%% file: paper_EMROM_basic.tex
\documentclass[a4paper]{article}

\input{packages.tex}
\usepackage[a4paper, total={6in, 8in}]{geometry}
\usepackage{csquotes}
\usepackage[backend=bibtex,style=numeric-comp,giveninits=true,doi=false,isbn=false,url=false,eprint=false,maxbibnames=99]{biblatex}

\input{macros.tex}

\graphicspath{{pictures/}}
\doublefalse

\title{{\papertitle}}
\input{header_basic.tex}

\begin{document}
	\maketitle

	\begin{abstract}
		\input{parts_abstract.tex}
	\end{abstract}


	\input{parts_intro.tex}
	\input{parts_methods.tex}
	\input{parts_results.tex}
	\input{parts_discussion.tex}
	\input{parts_conclusions.tex}

	\input{acknowledgements.tex}

	\begin{appendices}
		\input{parts_app_params.tex}
	\end{appendices}

	\printbibliography

\end{document}

%% file: packages.tex
\usepackage[english]{babel}
\usepackage[utf8]{inputenc}
\usepackage[colorinlistoftodos]{todonotes}
\usepackage{amsmath}
\usepackage{amssymb}
\usepackage{amsthm} 
\usepackage{mathtools}
\usepackage{mathrsfs}
\usepackage{siunitx}
\usepackage[hidelinks]{hyperref}
\usepackage{algorithm}
\usepackage{algpseudocode}

\usepackage{subfig}

\usepackage{tikz}

\usepackage{tabularx}
\usepackage{booktabs}
\usepackage{multirow}

\usepackage{empheq}

\usepackage{lipsum}

\usepackage[toc,page]{appendix}

\theoremstyle{plain}
\theoremstyle{plain}
\theoremstyle{plain}
\theoremstyle{plain}
\theoremstyle{plain}
\theoremstyle{definition}
\theoremstyle{remark}

\everymath{\displaystyle}

%% file: macros.tex
\newif\ifdouble

\newcommand{\MS}[1]{#1}

\newcommand{\papertitle}{A machine learning method for real-time numerical simulations of cardiac electromechanics}

\newcommand{\keywordOne}{Cardiac Electromechanics}
\newcommand{\keywordTwo}{Machine Learning}
\newcommand{\keywordThree}{Reduced Order Modeling}
\newcommand{\keywordFour}{Global Sensitivity Analysis}
\newcommand{\keywordFive}{Bayesian Parameter Estimation}

\newcommand{\expected}{\mathbb{E}}
\newcommand{\variance}{\mathbb{V}\mathrm{ar}}

\newcommand{\modEMfom}  {\mathcal{M}_{\text{3D}}}
\newcommand{\modEMred}  {\mathcal{M}_{\text{ANN}}}
\newcommand{\modCirc }  {\mathcal{C}}
\newcommand{\modEMCfom} {\modEMfom\textrm{-}\modCirc}
\newcommand{\modEMCred} {\modEMred\textrm{-}\modCirc}
\newcommand{\modCircW}  {\mathcal{C}'}
\newcommand{\modEMCWfom} {\modEMfom\textrm{-}\modCircW}
\newcommand{\modEMCWred} {\modEMred\textrm{-}\modCircW}

\newcommand{\modEMredA}  {\mathcal{M}_{\text{ANN}}^{\text{single}}}
\newcommand{\modEMredB}  {\mathcal{M}_{\text{ANN}}^{\text{full}}}
\newcommand{\modEMCredA} {\modEMredA\textrm{-}\modCirc}
\newcommand{\modEMCredB} {\modEMredB\textrm{-}\modCirc}
\newcommand{\modEMCWredA} {\modEMredA\textrm{-}\modCircW}
\newcommand{\modEMCWredB} {\modEMredB\textrm{-}\modCircW}

\newcommand{\THB}  {T_{\text{HB}}}

\newcommand{\EMState}  {\mathbf{y}}
\newcommand{\CircState}{\mathbf{c}}
\newcommand{\ANNState}{\mathbf{z}}
\newcommand{\EMRhs}  {\mathcal{L}}
\newcommand{\CircRhs}{\mathbf{f}}
\newcommand{\ANNRhs}{\mathcal{NN}}

\newcommand{\NumANNState} {N_{z}}
\newcommand{\NumANNWeights} {N_{w}}

\newcommand{\param}  {\mathbf{p}}
\newcommand{\paramC} {\param_{\mathcal{C}}}
\newcommand{\paramM} {\param_{\mathcal{M}}}

\newcommand{\paramSpace} {\mathscr{P}}
\newcommand{\paramSpaceC} {\paramSpace_{\mathcal{C}}}
\newcommand{\paramSpaceM} {\paramSpace_{\mathcal{M}}}
\newcommand{\NumParamC} {N_{\mathcal{C}}}
\newcommand{\NumParamM} {N_{\mathcal{M}}}

\newcommand{\error}{\boldsymbol{\epsilon}}
\newcommand{\qoi}  {\mathbf{q}}
\newcommand{\qoiObs}  {\qoi_{\text{obs}}}
\newcommand{\forward}  {\mathcal{F}}
\newcommand{\forwardRed}  {\widetilde{\mathcal{F}}}
\newcommand{\piPrior}{\pi_{\text{prior}}}
\newcommand{\piPost}{\pi_{\text{post}}}
\newcommand{\NoiseCov}{\boldsymbol{\Sigma}}

\newcommand{\errorROM}{\error_{\text{ROM}}}
\newcommand{\NoiseCovROM}{\NoiseCov_{\text{ROM}}}
\newcommand{\errorEXP}{\error_{\text{exp}}}
\newcommand{\NoiseCovEXP}{\NoiseCov_{\text{exp}}}
\newcommand{\NoiseMagnEXP}{\sigma_{\text{exp}}}

\newcommand{\VLVemFOM}{V_{\mathrm{LV}}^{\mathrm{3D}}}
\newcommand{\VLVemRED}{V_{\mathrm{LV}}^{\mathrm{ANN}}}
\newcommand{\VLVcirc} {V_{\mathrm{LV}}^{\mathrm{0D}}}

\newcommand{\prs}{p}
\newcommand{\vol}{V}

\newcommand{\PLVtrain}[1]{p_{\mathrm{LV}}^{#1}}
\newcommand{\VLVtrain}[1]{V_{\mathrm{LV}}^{#1}}

\newcommand{\paramCtrain}[1]{\param_{\mathcal{C}}^{#1}}
\newcommand{\paramMtrain}[1]{\param_{\mathcal{M}}^{#1}}

\newcommand{\numTrain}{N_{\text{train}}}
\newcommand{\Ttrain}{T_{\text{train}}}

\newcommand{\ANNparam}{\mathbf{w}}
\newcommand{\ANNparamTrained}{\widehat{\ANNparam}}
\newcommand{\ANNregularization}{\beta}
\newcommand{\ANNnumLayers}{N_{\text{layers}}}
\newcommand{\ANNnumNeurons}{N_{\text{neurons}}}

\newcommand{\EPcondf}{\sigma_{\mathrm{f}}}
\newcommand{\fibersAngle}{\alpha}
\newcommand{\MstiffPass}{C}

\newcommand{\Rtwo}{R\textsuperscript{2}}

\newcommand{\SobolFirst}[2]{S_{{#1}{#2}}}
\newcommand{\SobolTotal}[2]{S_{{#1}{#2}}^T}

\newcommand{\Vheart}{V_{\mathrm{heart}}^{\mathrm{tot}}}

\newcommand{\VLV}{V_{\mathrm{LV}}}

\newcommand{\VRV}{V_{\mathrm{RV}}}

\newcommand{\VnLA}{V_{\mathrm{0,LA}}}

\newcommand{\VnRA}{V_{\mathrm{0,RA}}}
\newcommand{\VnRV}{V_{\mathrm{0,RV}}}

\newcommand{\PLV}{p_{\mathrm{LV}}}

\newcommand{\PRV}{p_{\mathrm{RV}}}

\newcommand{\EpLA}{E_{\mathrm{LA}}^{\mathrm{pass}}}

\newcommand{\EpRA}{E_{\mathrm{RA}}^{\mathrm{pass}}}
\newcommand{\EpRV}{E_{\mathrm{RV}}^{\mathrm{pass}}}

\newcommand{\EaLA}{E_{\mathrm{LA}}^{\mathrm{act}}}

\newcommand{\EaRA}{E_{\mathrm{RA}}^{\mathrm{act}}}
\newcommand{\EaRV}{E_{\mathrm{RV}}^{\mathrm{act}}}

\newcommand{\tLAV}{t_{\mathrm{L}}^{\mathrm{av}}}
\newcommand{\tRAV}{t_{\mathrm{R}}^{\mathrm{av}}}
\newcommand{\TCLA}{T_{\mathrm{LA}}^{\mathrm{contr}}}

\newcommand{\TCRA}{T_{\mathrm{RA}}^{\mathrm{contr}}}
\newcommand{\TCRV}{T_{\mathrm{RV}}^{\mathrm{contr}}}
\newcommand{\TRLA}{T_{\mathrm{LA}}^{\mathrm{rel}}}

\newcommand{\TRRA}{T_{\mathrm{RA}}^{\mathrm{rel}}}
\newcommand{\TRRV}{T_{\mathrm{RV}}^{\mathrm{rel}}}

\newcommand{\CarSYS}{C_{\mathrm{AR}}^{\mathrm{SYS}}}
\newcommand{\CarPUL}{C_{\mathrm{AR}}^{\mathrm{PUL}}}
\newcommand{\CvnSYS}{C_{\mathrm{VEN}}^{\mathrm{SYS}}}
\newcommand{\CvnPUL}{C_{\mathrm{VEN}}^{\mathrm{PUL}}}

\newcommand{\RarSYS}{R_{\mathrm{AR}}^{\mathrm{SYS}}}
\newcommand{\RarPUL}{R_{\mathrm{AR}}^{\mathrm{PUL}}}
\newcommand{\RvnSYS}{R_{\mathrm{VEN}}^{\mathrm{SYS}}}
\newcommand{\RvnPUL}{R_{\mathrm{VEN}}^{\mathrm{PUL}}}

\newcommand{\LarSYS}{L_{\mathrm{AR}}^{\mathrm{SYS}}}
\newcommand{\LarPUL}{L_{\mathrm{AR}}^{\mathrm{PUL}}}
\newcommand{\LvnSYS}{L_{\mathrm{VEN}}^{\mathrm{SYS}}}
\newcommand{\LvnPUL}{L_{\mathrm{VEN}}^{\mathrm{PUL}}}

\newcommand{\Rmin}{R_{\mathrm{min}}}
\newcommand{\Rmax}{R_{\mathrm{max}}}

\newcommand{\aXB}{{a_{\text{XB}}}}

\newcommand{\VminLA}{\vol_{\mathrm{LA}}^{\min}}
\newcommand{\VmaxLA}{\vol_{\mathrm{LA}}^{\max}}
\newcommand{\pminLA}{\prs_{\mathrm{LA}}^{\min}}
\newcommand{\pmaxLA}{\prs_{\mathrm{LA}}^{\max}}
\newcommand{\VminLV}{\vol_{\mathrm{LV}}^{\min}}
\newcommand{\VmaxLV}{\vol_{\mathrm{LV}}^{\max}}
\newcommand{\pminLV}{\prs_{\mathrm{LV}}^{\min}}
\newcommand{\pmaxLV}{\prs_{\mathrm{LV}}^{\max}}
\newcommand{\SVLV}{\mathrm{SV}_{\mathrm{LV}}}
\newcommand{\VminRA}{\vol_{\mathrm{RA}}^{\min}}
\newcommand{\VmaxRA}{\vol_{\mathrm{RA}}^{\max}}
\newcommand{\pminRA}{\prs_{\mathrm{RA}}^{\min}}
\newcommand{\pmaxRA}{\prs_{\mathrm{RA}}^{\max}}
\newcommand{\VminRV}{\vol_{\mathrm{RV}}^{\min}}
\newcommand{\VmaxRV}{\vol_{\mathrm{RV}}^{\max}}
\newcommand{\pminRV}{\prs_{\mathrm{RV}}^{\min}}
\newcommand{\pmaxRV}{\prs_{\mathrm{RV}}^{\max}}
\newcommand{\SVRV}{\mathrm{SV}_{\mathrm{RV}}}
\newcommand{\pminARSYS}{\prs_{\mathrm{AR, SYS}}^{\min}}
\newcommand{\pmaxARSYS}{\prs_{\mathrm{AR, SYS}}^{\max}}


%% file: header_basic.tex
\author{Francesco Regazzoni$^{1,*}$,
		Matteo Salvador$^1$
		Luca Ded\`{e}$^1$,
		Alfio Quarteroni$^{1, 2}$}

\date{\footnotesize
	$^1$ MOX-Dipartimento di Matematica, Politecnico di Milano, Milan, Italy \\
	$^2$ \'Ecole Polytechnique F\'ed\'erale de Lausanne, Lausanne, Switzerland (\textit{Professor Emeritus})\\[2ex]
	$^*$ \textit{Corresponding author} (\texttt{francesco.regazzoni@polimi.it}) \\
    }

%% file: parts_abstract.tex
We propose a machine learning-based method to build a system of differential equations that approximates the dynamics of 3D electromechanical models for the human heart, accounting for the dependence on a set of parameters.
Specifically, our method permits to create a reduced-order model (ROM), written as a system of Ordinary Differential Equations (ODEs) wherein the forcing term, given by the right-hand side, consists of an Artificial Neural Network (ANN), that possibly depends on a set of parameters associated with the electromechanical model to be surrogated.
This method is non-intrusive, as it only requires a collection of pressure and volume transients obtained from the full-order model (FOM) of cardiac electromechanics.
Once trained, the ANN-based ROM can be coupled with hemodynamic models for the blood circulation external to the heart, in the same manner as the original electromechanical model, but at a dramatically lower computational cost.
Indeed, our method allows for real-time numerical simulations of the cardiac function.
Our results show that the ANN-based ROM is accurate with respect to the FOM (relative error between $10^{-3}$ and $10^{-2}$ for biomarkers of clinical interest), while requiring very small training datasets (30-40 samples).
We demonstrate the effectiveness of the proposed method on two relevant contexts in cardiac modeling.
First, we employ the ANN-based ROM to perform a global sensitivity analysis on both the electromechanical and hemodynamic models.
Second, we perform a Bayesian estimation of two parameters starting from noisy measurements of two scalar outputs.
In both these cases, replacing the FOM of cardiac electromechanics with the ANN-based ROM makes it possible to perform in a few hours of computational time all the numerical simulations that would be otherwise unaffordable, because of their overwhelming computational cost, if carried out with the FOM.
As a matter of fact, our ANN-based ROM is able to speedup the numerical simulations by more than three orders of magnitude.

%% file: parts_intro.tex
\section{Introduction}
\label{sec:introduction}

Numerical simulations of cardiac electromechanics are gaining momentum in the context of cardiovascular research and computational medicine \cite{smith2004multiscale,vigmond2008effect,trayanova2011whole,fink2011cardiac,washio2015ventricular,quarteroni2019cardiovascular,dede2021mathematical}.
Cardiac in silico models are based on an anatomically detailed and biophysically accurate representation of the human heart, consisting of multiscale mathematical models of several physical processes, from the cellular to the organ scale, described by systems of differential equations.
However, the clinical exploitation of cardiac numerical simulations is seriously hampered by their overwhelming computational cost.
As a matter of fact, the simulation of a single heartbeat for an anatomically accurate patient-specific model may require several hours of computational time even on a supercomputer platform.

A promising approach to address this issue is to replace the computationally expensive cardiac \MS{electromechanical} model, say the full-order model (FOM), with a reduced version of it, called reduced-order model (ROM), to be called any time new parameters come in.
In the so-called \textit{offline phase}, the ROM is built from a database of numerical simulations that are previously obtained by solving the FOM itself.
Then, in the \textit{online phase}, the ROM is used as a surrogate of the FOM, at a dramatically reduced computational cost, for any new instance of the parameters.
Clearly, this procedure pays off if the computational gain of the online phase outweighs the cost of the offline phase, or whenever the online phase requires real-time execution, which is often the case in the clinical practice, where timeliness is a key factor.

Recently, this framework has been applied in the context of cardiac modeling, primarily by using machine learning algorithms, including Gaussian Process emulators (GPEs), Artificial Neural Networks (ANNs), decision tree algorithms such as eXtreme Gradient Boosting (XGBoost) and K-Nearest Neighbor (KNN) \cite{di2018gaussian,dabiri2019prediction,Longobardi2020,Cai2021,coveney2021bayesian}.
These emulators are trained to fit the map that links the model parameters with a set of scalar outputs of interest, known as quantities of interest (QoIs), which represent clinically meaningful biomarkers.
This map is learned from a collection of numerical simulations, obtained by sampling the parameter space, which can then be used as FOM surrogate.

Replacing a FOM with a ROM is particularly advantageous in the so-called many-query settings, i.e. when the forward model has to be solved many times for different parameter values.
This is the case, for example, of parameter calibration \cite{brooks1998markov,hoffman2013stochastic} and sensitivity analysis \cite{plischke2013global,sobol1990sensitivity,xu2008uncertainty,kucherenko2012estimation}, fundamental tools for a reliable use of computational models in the clinical practice.
In fact, only a few of the many parameters feeding cardiac electromechanical models can be directly measured, while most of them must be estimated through indirect (possibly non-invasive) measurements, by solving suitable (typically ill-posed) inverse problems \cite{aster2018parameter}.
It is therefore important to quantify how much each parameter affects the outputs of the model and, vice versa, how uncertainty on measured data \MS{reverberates on parameter uncertainty}.
These requirements are at the basis of the \textit{verification and validation under uncertainty} (VVUQ) paradigm \cite{oberkampf2004verification,hu20162014} and are often necessary to comply with standards issued by health agencies and policy makers, such as the ASME V\&V40 standard \cite{asme2018assessing} recognized by the FDA \cite{fda}.


In this paper, we propose a machine learning \MS{method} to build ROMs of cardiac \MS{electromechanical} models.
Our approach relies on the ANN-based method that we proposed in \cite{regazzoni2019modellearning}, which can learn a time-dependent differential equation from a collection of input-output pairs.
With respect to existing approaches, we only surrogate the time-dependent pressure-volume relationship of a cardiac chamber, while we do not reduce the model describing external circulation.
The latter is indeed either a low dimensional 0D windkessel or closed-loop circulation model comprised of few state variables (up to two dozens), which does not require further reduction.
In other terms, we derive an ANN-based ROM for the computationally demanding 3D components, while we leave in FOM version the lightweight ones.
Unlike emulators, for which the online phase consists in evaluations of the map linking model parameters to QoIs, with our approach the online phase consists instead in numerical simulations, in which the ANN-based ROM of the electromechanical model is coupled with the circulation model, at a very low computational cost.
As a matter of fact, these numerical simulations can be performed in real-time on a standard laptop.

We consider two different test cases to prove the efficacy of our mathematical approach.
We perform variance-based global sensitivity analysis on both electromechanical and hemodynamic parameters and Bayesian estimation by means of the Markov Chain Monte Carlo (MCMC) method to infer two parameters from noisy measurements of two scalar outputs.

This paper is organized as follows.
In Sec.~\ref{sec:methods} we present the proposed method and the models used to produce the numerical results, which are presented and commented in Sec.~\ref{sec:results}.
Then, in Sec.~\ref{sec:discussion}, we critically discuss the obtained results and the pros and cons of the proposed method, with respect to other methods available in the literature.
Finally, we draw our conclusions and final remarks in Sec.~\ref{sec:conclusions}.

This manuscript is accompanied by \textcolor{blue}{\url{https://github.com/FrancescoRegazzoni/cardioEM-learning}}, a public repository containing the codes and the datasets necessary to reproduce the presented results.

%% file: parts_methods.tex
\section{Methods}
\label{sec:methods}

We introduce our method in an abstract formulation for cardiac electromechanics.
As a matter of fact, our approach, thanks to its non-intrusive nature, can be applied to different \MS{electromechanical} models that are available in literature.

\subsection{The full-order model (FOM)}
\label{sec:methods:FOM}

Let us consider a generic model of cardiac electromechanics, that is a set of differential equations describing physical processes involved in the heart function.
We introduce the state vector $\EMState(t)$, collecting the state variables associated with this multiphysics system.
These may include the transmembrane potential, gating variables, ionic concentrations, protein states, tissue displacement, or simply phenomenological variables.
In this paper, we focus on the single-chamber case of the human heart (e.g. the left ventricle, that we call from now on LV), as the generalization to multiple chambers is straightforward.
By introducing a nonlinear differential operator $\EMRhs$ that encodes the differential equations and boundary conditions associated with the electromechanical model, the latter reads
\begin{equation} \label{eqn:model_3D}
    \left\{
    \begin{aligned}
        \frac{\partial \EMState(t)}{\partial t} &= \EMRhs  (\EMState(t)  , \PLV(t), t; \paramM) && \text{for } t \in (0,T],\\
        \EMState  (0) &= \EMState_0,   && \\
    \end{aligned}
    \right.
\end{equation}
where $\PLV(t)$ denotes the LV endocardial pressure (here seen as an input), $\paramM$ are the model parameters (possibly including, e.g., electrical conductivities, cell membrane conductances, protein binding affinities, contractility, passive tissue properties) and $\EMState_0$ is the initial state.
We denote by $\paramSpaceM \subseteq \mathbb{R}^{\NumParamM}$ the space of parameters such that $\paramM \in \paramSpaceM$, being $\NumParamM$ the number of parameters.
We notice that the right-hand side of \eqref{eqn:model_3D} depends on $t$, as heartbeats are paced by externally applied stimuli that we assume to have a period of duration $\THB$, which is fixed a priori.

\begin{figure}
	\centering
	\includegraphics[width=\textwidth]{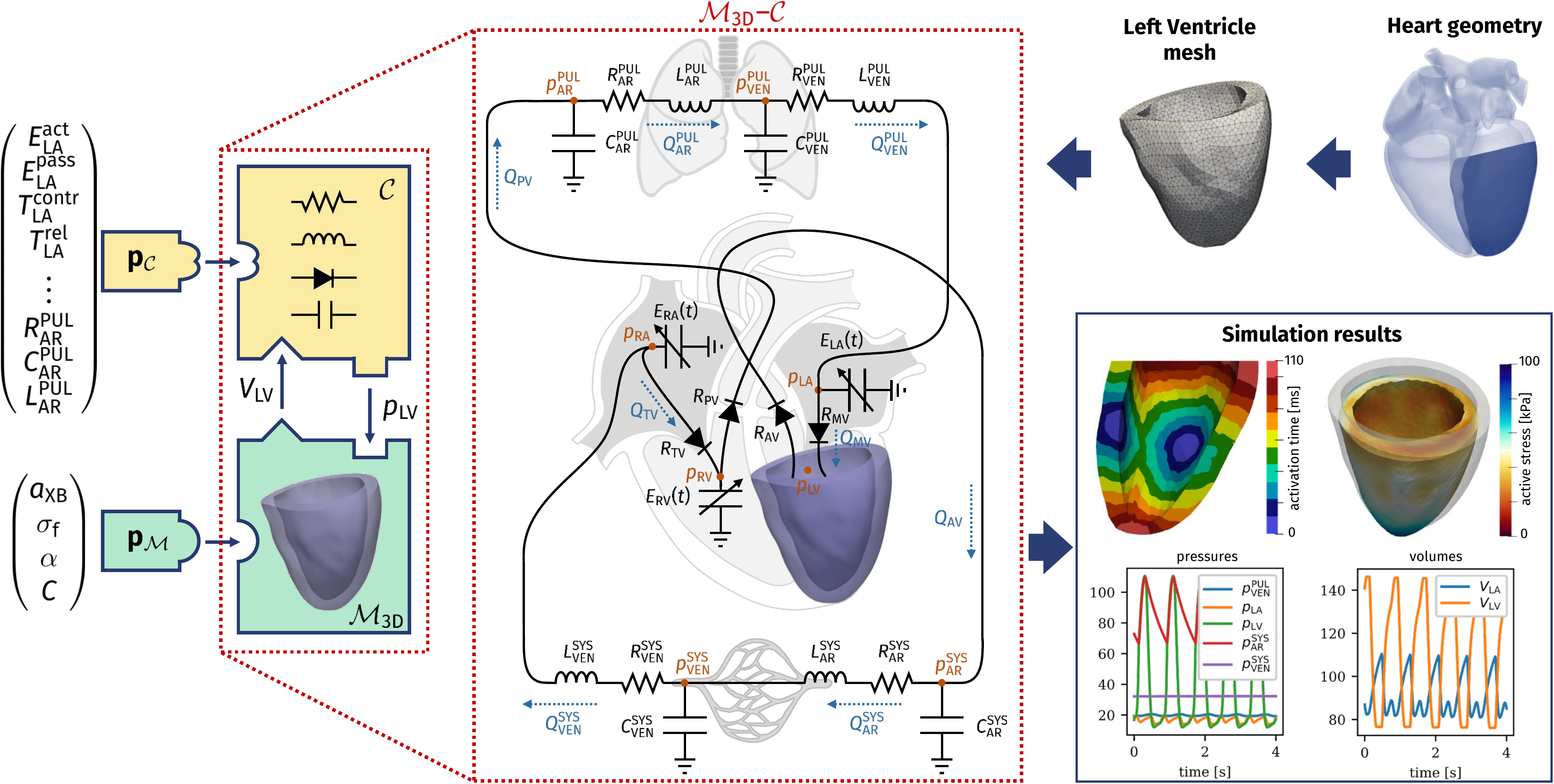}
	\caption{
    The $\modEMCfom$ model.
    The parameters $\paramC$ and $\paramM$ (left) are associated with the $\modCirc$ and $\modEMfom$ respectively.
    The two mathematical models are coupled via the variables $\PLV$ and $\VLV$.
    Their union constitutes the model $\modEMCfom$.
    The model $\modEMCfom$ that we consider to produce the numerical results of this paper is shown in the center of the figure. 
    For more details on the $\modEMCfom$ model and for the definition of the parameters $\paramC$ and $\paramM$, see Sec.~\ref{sec:methods:EM_model}.
    The 3D Finite Element model requires a computational mesh, obtained from the patient's heart geometry (top right).
    The output of the numerical simulations (bottom right) consists in spatially distributed fields (such as activation time, active tension and tissue displacement) and in 0D transients (such as pressures and volumes).
    }
	\label{fig:model_MC}
\end{figure}

The 3D cardiac \MS{electromechanical} model \eqref{eqn:model_3D}, henceforth denoted by $\modEMfom$, must be coupled with a closure relationship assigning the pressure $\PLV(t)$.
One possible option is to couple the $\modEMfom$ model with a 0D model for the external circulation (see e.g. \cite{hirschvogel2017monolithic,regazzoni2021em_circulation_pt1,augustin2020physiologically}), thus obtaining a system in the form of
\begin{equation} \label{eqn:model_3D-0D}
    \left\{
    \begin{aligned}
        \frac{\partial \EMState(t)}{\partial t} &= \EMRhs  (\EMState(t)  , \PLV(t), t; \paramM) && \text{for } t \in (0,T],\\
        \frac{d \CircState(t)}{dt} &= \CircRhs(\CircState(t), \PLV(t), t; \paramC) && \text{for } t \in (0,T],\\
        \VLVcirc(\CircState(t)) &= \VLVemFOM(\EMState(t))         && \text{for } t \in (0,T],\\
        \EMState  (0) &= \EMState_0,   && \\
        \CircState(0) &= \CircState_0, && \\
    \end{aligned}
    \right.
\end{equation}
where $\CircState(t)$ are the state variables of the circulation model (pressures, volumes and fluxes in the circulatory network) and $\paramC \in \paramSpaceC \subseteq \mathbb{R}^{\NumParamC}$ is a vector of $\NumParamC$ parameters (e.g. vascular resistances and conductances).
The mechanical and hemodynamic models are coupled through the geometric consistency relationship $\VLVemFOM(\EMState(t)) = \VLVcirc(\CircState(t))$, where the left-hand and right-hand sides represent the LV volume predicted by the $\modEMfom$ and by the $\modCirc$ models, respectively.
The LV pressure $\PLV$ is determined as a Lagrange multiplier that enforces the consistency relationship.
An alternative approach is to adopt different closure relationships in the different phases of the heartbeat \cite{levrero2020sensitivity, Gerbi2018monolithic} with suitable preload and afterload models, such as windkessel models \cite{westerhof2009arterial}.
These relationships link the changes in LV pressure $\PLV$ with its volume, obtained as $\VLV = \VLVemFOM(\EMState)$.
In both cases, should the $\modEMfom$ be coupled to a closed-loop circulation model or to an afterload-preload relationship, we denote by $\modEMCfom$ the resulting coupled model.
In Fig.~\ref{fig:model_MC} we show an $\modEMCfom$ model in which $\modCirc$ consists of a the lumped-parameter  closed-loop model of \cite{regazzoni2021em_circulation_pt1}, which is employed to produce the numerical results of Sec.~\ref{sec:results}.

We remark that the closure relationships never directly involve the state $\EMState$ of the model $\modEMfom$, but only the LV volume $\VLV = \VLVemFOM(\EMState)$.
This is a key observation since it suggests that a ROM that is able to surrogate the relationship between $\PLV$ and $\VLV$, albeit agnostic of the state $\EMState$, can replace the $\modEMfom$ model in its coupling to the $\modCirc$ model.
In what follows, we present our strategy to build an ANN-based ROM of the $\modEMfom$ model, denoted by $\modEMred$, which can be coupled with the $\modCirc$ model resulting in the $\modEMCred$ coupled model, that will surrogate the $\modEMCfom$ model.

\subsection{The reduced-order model (ROM)}
\label{sec:methods:ROM}

To setup a ROM surrogating the $\modEMfom$ model of Eq.~\eqref{eqn:model_3D}, we employ the machine learning method that we proposed in \cite{regazzoni2019modellearning}.
This method is designed to learn a differential equation from time-dependent input-output pairs, by training an Ordinary Differential Equation (ODE) model, whose right-hand side is represented by an ANN.
In this work, we define our ANN-based ROM $\modEMred$ as
\begin{equation} \label{eqn:model_ANN}
    \left\{
    \begin{aligned}
        \frac{d \ANNState(t)}{d t} &= \ANNRhs\left(
            \ANNState(t),
            \PLV(t),
            \cos(\tfrac{2 \pi t}{\THB}),
            \sin(\tfrac{2 \pi t}{\THB}),
            \paramM;
            \ANNparamTrained
        \right)
        && \text{for } t \in (0,T],\\
        \ANNState  (0) &= \ANNState_0,   && \\
    \end{aligned}
    \right.
\end{equation}
where $\ANNState(t) \in \mathbb{R}^{\NumANNState}$ is the reduced state and $\ANNRhs \colon \mathbb{R}^{\NumANNState + \NumParamM + 3} \to \mathbb{R}^{\NumANNState}$ is a fully connected ANN.
The ANN input consists indeed of $\NumANNState$ state variables, $\NumParamM$ scalar parameters, the pressure $\PLV$, and the two periodic inputs $\cos({2 \pi t}/{\THB})$ and $\sin({2 \pi t}/{\THB})$ (whose role will be clarified later), for a total of $\NumANNState + \NumParamM + 3$ input neurons.
The vector $\ANNparamTrained \in \mathbb{R}^{\NumANNWeights}$ \MS{encodes} the weights and biases of the ANN, that need to be suitably trained.
We remark that, among the arguments of the ANN, we have not introduced the time variable $t$, but rather $\cos({2 \pi t}/{\THB})$ and $\sin({2 \pi t}/{\THB})$, that are the coordinates of a point cyclically moving along a circumference with period $\THB$.
In this way, the ROM encodes by construction the periodicity associated with the heartbeat pacing.
This expedient allows for the use of the ROM also for time spans longer than those shown during the training phase.
Moreover, by introducing the parametric dependence (i.e. on $\paramM$) within the ANN, the latter is not specific to a particular parameter setting.

In this work, according to \cite{regazzoni2019modellearning}, we adopt an output-inside-the-state approach, that is we train the model so that the LV volume coincides with the first state variable.
More precisely, the LV volume predicted by the $\modEMred$ model is by definition $\VLVemRED(\ANNState(t)) := \ANNState(t)  \cdot \mathbf{e}_1$, where $\mathbf{e}_1$ is the first element of the canonical basis of $\mathbb{R}^{\NumANNState}$.
Therefore, the first ROM state variable has a clear physical interpretation (i.e. it coincides with the LV volume), while the other ROM states are latent variables with no immediate physical interpretation, providing \MS{however} a compact representation of the full-order state $\EMState(t)$.
Coherently with this choice, we define the initial state as $\ANNState_0 = (\VLVemFOM(\EMState_0), \mathbf{0})^T$.
The remaining initial states are set, without loss of generality, to zero (this choice does not reduce the space of candidate models, as proved in \cite{regazzoni2019modellearning}).

\begin{figure}
	\centering
	\includegraphics[width=0.9\textwidth]{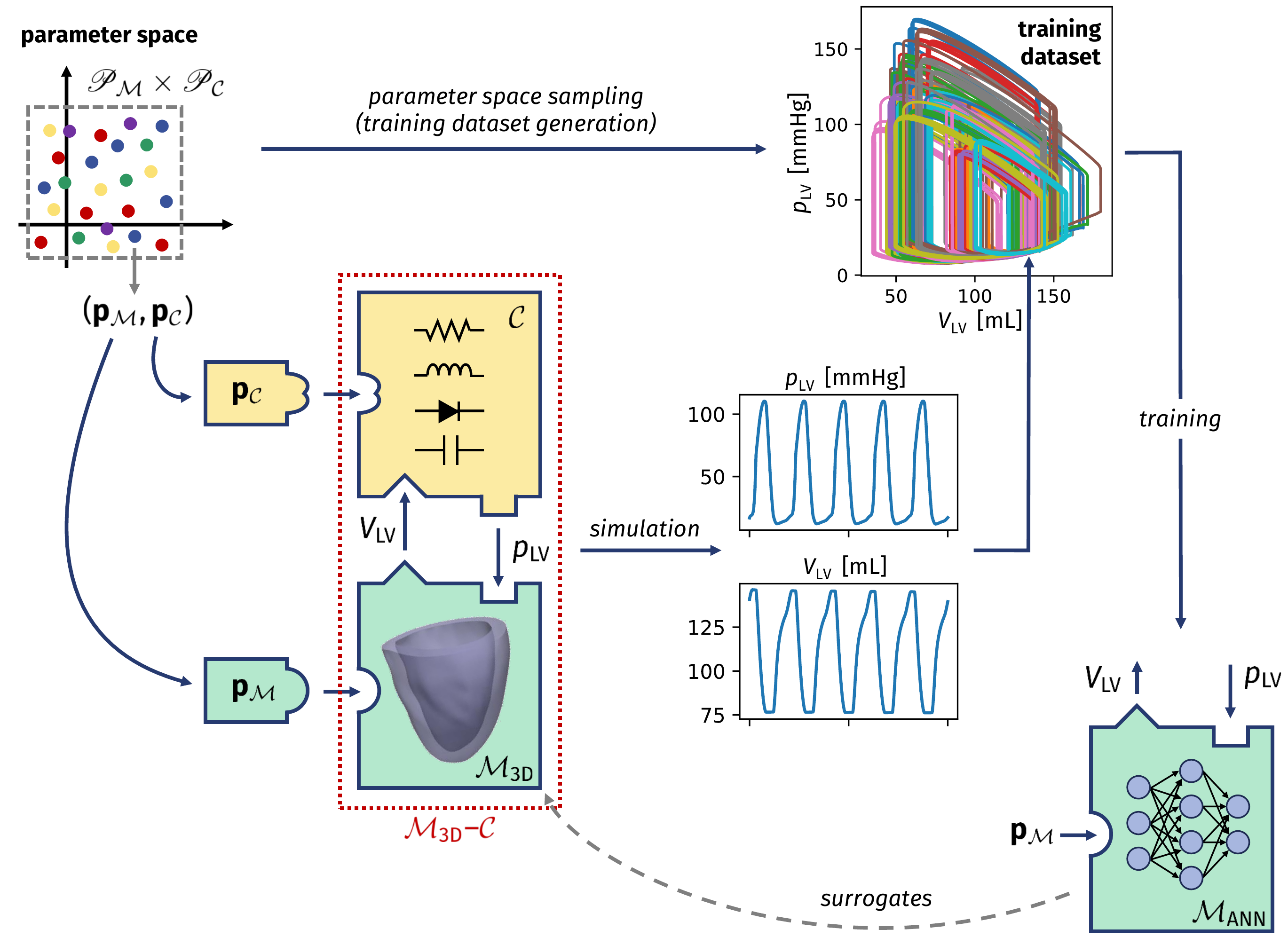}
	\caption{
	Training pipeline of the $\modEMred$ model.
	First, we sample the parameter space $\paramSpaceM \times \paramSpaceC$ (top left figure) and, for each parameter instance $(\paramM, \paramC)$, we simulate some heartbeats through the $\modEMCfom$ model (center figure).
	Finally, from the training set obtained by collecting the resulting pressure and volume transients (top right figure), we train the ANN-based model $\modEMred$ (bottom right figure), according to Eq.~\eqref{eqn:optimization_problem}.
	}
	\label{fig:training_flow}
\end{figure}

Here we need to determine the optimal value of the weights $\ANNparamTrained$, such that the $\modEMred$ model reproduces the outputs of the $\modEMfom$ model as accurately as possible.
With this goal, we generate a training set, by sampling the parameter space $\paramSpaceM \times \paramSpaceC$ with $\numTrain$ sample points.
For each sample, we perform a simulation with the $\modEMCfom$ model until time $\Ttrain$, and we record the LV pressure and volume transients.
The training set is thus given by
$$
\paramCtrain{i}, \paramMtrain{i}, \PLVtrain{i}(t), \VLVtrain{i}(t)
\qquad
t \in [0, \Ttrain], \qquad \text{for } i = 1, \dots, \numTrain.
$$
We remark that, due to the non-intrusive nature of our method, there is no need to retain the FOM states $\EMState(t)$.
Finally, we train the ANN weights $\ANNparamTrained$ by minimizing the discrepancy between the training data and the model outputs, that is by considering the following constrained optimization problem
\begin{equation} \label{eqn:optimization_problem}
    \left\{
    \begin{aligned}
        &\ANNparamTrained =
        \underset{\ANNparam \in \mathbb{R}^{\NumANNWeights}}{\operatorname{argmin}}
        \left[
            \sum_{i = 1}^{\numTrain} \int_0^{\Ttrain} | \VLVtrain{i}(t) - \VLVemRED(\ANNState^i(t)) |^2 \, dt
            + \ANNregularization |\ANNparam|^2
        \right]
        \\
        &\text{such that, for each $i = 1, \dots, \numTrain$:}
        \\[1mm] &\qquad
        \begin{aligned}
        \frac{d \ANNState^i(t)}{dt} &= \ANNRhs\left(
            \ANNState^i(t),
            \PLVtrain{i}(t),
            \cos(\tfrac{2 \pi t}{\THB}),
            \sin(\tfrac{2 \pi t}{\THB});
            \paramMtrain{i};
            \ANNparam
        \right)
        &&
        \text{for } t \in (0, \Ttrain],
        \\
        \ANNState^i(0) &= \ANNState_0,
        &&
        \end{aligned}
    \end{aligned}
    \right.
\end{equation}
where $\ANNregularization > 0$ is a regularization hyperparameter.
We remark that the optimization problem \eqref{eqn:optimization_problem} is not a standard machine learning problem of data fitting.
Indeed, the ANN appears at the right-hand side of a differential equation that acts as a constraint under which the loss function is minimized.
To train this model, we use the algorithm proposed in \cite{regazzoni2019modellearning}, which envisages approximating the differential equation by the Forward Euler method, the loss function by the trapezoidal method, and then computing the gradients by solving the adjoint equations.
The parameters are then optimized by means of the Levenberg-Marquardt method \cite{nocedal2006numopt}.
The training pipeline is summarized in Fig.~\ref{fig:training_flow}.

\begin{figure}
	\centering
	\includegraphics[width=\textwidth]{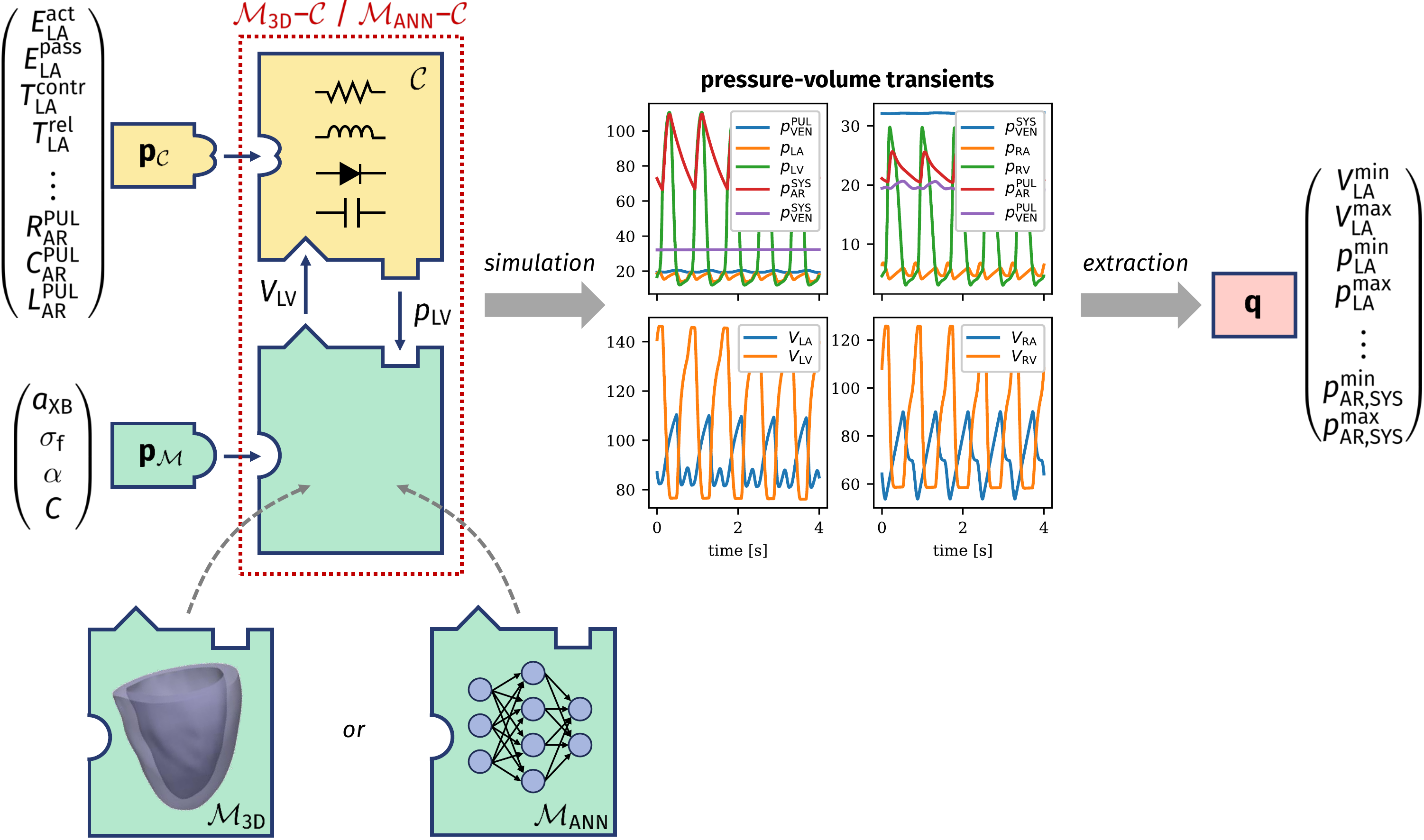}
	\caption{
	Parameters-to-QoIs computation using either the $\modEMfom$ or the $\modEMred$ model.
    Given a parameter instance $(\paramM, \paramC)$, either the $\modEMfom$ or the $\modEMred$ model can be coupled with the $\modCirc$ model to obtain pressure and volume transients, from which a set of QoIs are extracted.
    See Tabs.~\ref{tab:paramsM}, \ref{tab:paramsC} and \ref{tab:QoI} for the definition of $\paramM$, $\paramC$ and $\qoi$, respectively.}
	\label{fig:param_qoi_flow}
\end{figure}

Once the ANN has been trained (that is, the optimal weights $\ANNparamTrained$ have been determined), the $\modEMred$ model can be used as a surrogate of the $\modEMfom$ model, also for different combinations of the parameters than those contained in the training set.
Moreover, it can be coupled with the closure relationships $\modCirc$, thus obtaining the $\modEMCred$ model.
For example, by considering the case of a closed-loop circulation model, as in \eqref{eqn:model_3D-0D}, its reduced counterpart $\modEMCred$ reads
\begin{equation} \label{eqn:model_ANN-0D}
    \left\{
    \begin{aligned}
        \frac{d \ANNState(t)}{d t} &= \ANNRhs\left(
            \ANNState(t),
            \PLV(t),
            \cos(\tfrac{2 \pi t}{\THB}),
            \sin(\tfrac{2 \pi t}{\THB});
            \paramM;
            \ANNparamTrained
        \right)
        && \text{for } t \in (0,T],\\
        \frac{d \CircState(t)}{dt} &= \CircRhs(\CircState(t), \PLV(t), t; \paramC)
        && \text{for } t \in (0,T],
        \\
        \VLVcirc(\CircState(t)) &= \VLVemRED(\ANNState(t))
        &&
        \text{for } t \in (0,T],
        \\
        \ANNState  (0) &= \ANNState_0,
        && \\
        \CircState(0) &= \CircState_0.
        && \\
    \end{aligned}
    \right.
\end{equation}
As matter of fact, the $\modEMfom$ model has the same input-output structure of the model being surrogated, that is $\modEMred$.
Hence, the latter can be employed in replacement of the former in approximating the outputs associated with a given set of parameters $(\paramM, \paramC)$, as shown in Fig.~\ref{fig:param_qoi_flow}.

\subsection{Hyperparameters tuning}
\label{sec:methods:hyperparameters}

Like any machine learning algorithm, our method depends on a set of \textit{hyperparameters}, that is variables that are not trained (they are not part of the vector $\ANNparamTrained$), but are used to control the training process.
These include ANN architecture hyperparameters (namely the number of layers $\ANNnumLayers$ and the number of neurons per layer $\ANNnumNeurons$), the regularization factor $\ANNregularization$ and the number of reduced states $\NumANNState$.
To tune the hyperparameters, we rely on a $k$-fold cross-validation procedure.
Specifically, after splitting the training set into $k = 5$ non-overlapping subsets, we cyclically train the model by excluding one subset that is used as validation set.
Finally, we evaluate the average validation errors and we select the hyperparameters setting that attains the lowest validation error and better generalization properties (see Appendix~\ref{app:params:Xvalidation} for more details).

\subsection{Global sensitivity analysis}
\label{sec:methods:GSA}

\MS{To assess how much each parameter contributes to the determination of a QoI, e.g. a biomarker of clinical interest, we perform a global sensitivity analysis.
This is typically done by sampling the parameter space and by computing suitable indicators, such as Borgonovo indices \cite{plischke2013global}, Sobol indices \cite{sobol1990sensitivity,homma1996importance}, Morris elementary effects \cite{morris1991factorial}, ANCOVA indices \cite{xu2008uncertainty} and Kucherenko indices \cite{kucherenko2012estimation}.
In this work, we perform a variance-based sensitivity analysis, which relies on a probabilistic approach.
We compute Sobol indices to quantify the sensitivity of a QoI (say $\qoi_j$) with respect to a parameter (say $\param_i$).}
Specifically, the so-called \textit{first-order Sobol index} (denoted by $\SobolFirst{i}{j}$) indicates the impact on the $j$-th QoI (i.e. $\qoi_j$) of the $i$-th parameter (i.e. $\param_i$) when the latter varies alone, according to the definition
\begin{equation*}
\SobolFirst{i}{j} = \frac{\variance_{\param_i}\left[ \expected_{\param_{\sim i}}\left[ \qoi_j | \param_i\right] \right]}{\variance\left[ \qoi_j \right]},
\end{equation*}
where $\param_{\sim i}$ indicates the set of all parameters excluding the $i$-th one.
The first-order Sobol index $\SobolFirst{i}{j}$, however, only accounts for the variations of the parameter $\param_i$ alone, averaged over variations in the other parameters, and thus does not account for the \textit{interactions} among parameters.
Conversely, it is possible to assess the importance of a parameter in determining a QoI, also accounting for the interactions among parameters, by resorting to the so-called \textit{total-effect Sobol index} $\SobolTotal{i}{j}$, defined as
\begin{equation*}
    \SobolTotal{i}{j} = \frac{\expected_{\param_{\sim i}}\left[ \variance_{\param_i}\left[ \qoi_j | \param_{\sim i}\right] \right]}{\variance\left[ \qoi_j \right]}
                  = 1 - \frac{\variance_{\param_{\sim i}}\left[ \expected_{\param_i}\left[ \qoi_j | \param_{\sim i}\right] \right]}{\variance\left[ \qoi_j \right]}.
\end{equation*}
The latter index quantifies the impact of a given parameter when it varies alone or together with other parameters \cite{sobol1990sensitivity}.

To compute an estimate of the Sobol indices, we use the Saltelli method \cite{homma1996importance,saltelli2002making}, that makes use of Sobol quasi-random sequences to approximate the integrals that need to be computed.
In practice, this method requires evaluating the model for a large number of parameters, and then processing the obtained QoIs to provide an estimate of the Sobol indices.

In this paper we perform a variance-based global sensitivity analysis simultaneously with respect to the circulation model parameters $\paramC$ and the electromechanical model parameters $\paramC$ (that is, we set $\param = (\paramM, \paramC)$).
To this end, we use model $\modEMCred$ as a surrogate for model $\modEMCfom$ to perform the evaluations required by the Saltelli method.
As we will see in Sec.~\ref{sec:results}, using $\modEMCred$ instead of $\modEMCfom$ entails a huge computational gain.
In addition, for each parameter setting we simulate a certain number of heartbeats to achieve a limit cycle (i.e. a periodic solution), and we calculate QoIs with respect to the last cycle, which is the most significant one (as it removes the effects of incorrect initializations of the state variables).
We remark that the need to simulate a certain number of heartbeats to overpass the transient phase (that typically lasts from 5 to 15 cycles, see Appendix~\ref{app:params:limit-cycle}) makes the use of a reduced cardiac \MS{electromechanical} model even more necessary.

\subsection{Parameter estimation under uncertainty}
\label{sec:methods:bayesianPE}


The patient-specific personalization of a cardiac \MS{electromechanical} model requires, besides the usage of a geometry derived from imaging data, the estimation of the parameters associated with the mathematical model (or at least the most important ones), starting from clinical measurements.
Very often, only a few scalar quantities are available for this purpose; moreover, the resolution of this inverse problem (i.e. estimating $\param$ from $\qoi$) should account for the noise that unavoidably affects the measurement of $\qoi$ and that reflects in uncertainty on $\param$.

Bayesian methods, such as MCMC \cite{brooks1998markov} and Variational Inference \cite{hoffman2013stochastic}, permit to address these issues within a rigorous statistical framework, by providing the \textit{likelihood}, expressed as a probability distribution of the parameters values, given the observed QoIs (denoted by $\qoiObs$).
Unlike methods that provide point estimates, either based on gradient-descent \cite{nocedal2006numopt} or genetic \cite{holland1992genetic} optimization algorithms, Bayesian methods are able to provide a probability distribution on the parameter space, encoding the \textit{credibility} of each parameter combination.
The credibility computation accounts for the uncertainty on measurements (associated with measurement noise and encoded in the noise \MS{covariance} matrix $\NoiseCov$), as well as a \textit{prior distribution} on parameters (denoted by $\piPrior(\param)$), that is previous knowledge or belief about the parameters.
By denoting the parameters-to-QoIs map by $\forward \colon \param \mapsto \qoi$, we have $\qoiObs = \forward(\param) + \error$, where $\error \sim \mathcal{N}( \cdot | \mathbf{0}, \NoiseCov)$ denotes the measurement error (that we assume for simplicity to be distributed as a Gaussian random variable).
Bayes' theorem \MS{states} that the \textit{posterior distribution} of parameters, that is the degree of belief of their value after having observed $\qoiObs$, is given by
\begin{equation*}
    \piPost(\param) = \frac{1}{Z}\, \mathcal{N}( \qoiObs | \forward(\param), \NoiseCov) \, \piPrior(\param),
\end{equation*}
where the normalization constant $Z$ is defined as
\begin{equation*}
    Z = \int_{\paramSpace} \mathcal{N}( \qoiObs | \forward(\widehat{\param}), \NoiseCov) \, d\piPrior(\widehat{\param}).
\end{equation*}
In practice, the computation of $\piPost$ may be challenging from the computational viewpoint, because of the need to approximate the integral that defines $Z$.
The \MS{MCMC} method permits to approximate the distribution $\piPost$ with a relatively small computational effort.
Similarly to the Saltelli method that we use for sensitivity analysis, also the MCMC method requires a large number of model evaluations, for different parameter values.
Moreover, this method is non intrusive, that is it only requires evaluations of the map $\forward \colon \param \mapsto \qoi$.
Therefore, we can employ for this purpose the $\modEMCred$ model as a surrogate for the $\modEMCfom$ model, which drastically reduces the necessary computational time.

Moreover, we remark that the Bayesian framework permits to rigorously account for the approximation error introduced by replacing the high fidelity model $\modEMCfom$ by its surrogate $\modEMCred$.
Indeed, if we denote by $\forwardRed$ the approximated parameters-to-QoIs map represented by the surrogate model $\modEMCred$, we have $\forward(\param) = \forwardRed(\param) + \errorROM$, where $\errorROM$ is the ROM approximation error.
It follows $\qoiObs = \forwardRed(\param) + \errorROM + \errorEXP$, where $\errorEXP$ is the experimental measurement error.
Since the two sources of error can be assumed independent, the covariance of the total error $\error = \errorROM + \errorEXP$ satisfies $\NoiseCov = \NoiseCovROM + \NoiseCovEXP$, where $\NoiseCovROM$ is the ROM approximation error covariance (which can be estimated by evaluating the ROM on a testing set) and $\NoiseCovEXP$ is the experimental error covariance (that depends on the measurement protocol at hand).
This permits to take into account in the estimation process the error introduced by the surrogate model.

To \MS{assess the capability} of our ROM to accelerate the estimation of parameters for multiscale cardiac \MS{electromechanical} models, we perform the following test.
First, we perform a simulation with the $\modEMCfom$ model, from which we derive a set of QoIs ($\qoiObs$).
Then, by employing the $\modEMCred$ model as a surrogate of the $\modEMCfom$ model, we obtain a Bayesian estimate of the parameters, that we validate against the values used to generate $\qoiObs$.

\subsection{The cardiac \MS{electromechanical} model}
\label{sec:methods:EM_model}

As mentioned above, due to its non-intrusive nature, our method is not limited to a specific cardiac \MS{electromechanical} model, but can be applied to virtually any electromechanical model as long as pressure and volume transients are available for the training procedure.
In this section, we briefly introduce the specific model used to produce the numerical results of this paper.

We consider a LV geometry processed from the Zygote 3D heart model \cite{zygote} endowed with a fiber architecture generated by means of the Bayer-Blake-Plank-Trayanova algorithm \cite{bayer2012novel,piersanti2020modeling}.
Before starting the simulations, we recover the reference unstressed configuration through the algorithm that we proposed in \cite{regazzoni2021em_circulation_pt1}.
To model the action potential propagation, we employ the monodomain equation \cite{collifranzone2014book}, coupled with the \textit{ten Tusscher-Panfilov} ionic model \cite{ten2006alternans}.
We model the microscale generation of active force through the biophysically detailed RDQ20-MF model \cite{regazzoni2020biophysically}, that is coupled, within an active stress approach, with the elastodynamics equations describing tissue mechanics.
On the other hand, the passive behavior of the tissue is modeled through a quasi-incompressible exponential constitutive law \cite{usyk2002computational}.
We model the interaction with the pericardium by means of spring-damper boundary conditions at the epicardium, while we adopt energy-consistent boundary conditions \cite{regazzoni2019mor-sarcomeres} to model the interaction with the part of the myocardium beyond the artificial ventricular base.
To model blood circulation, that is $\modCirc$, we rely on the 0D closed-loop model presented in \cite{regazzoni2021em_circulation_pt1}, consisting of a compartmental description of the cardiac chambers, systemic and pulmonary, arterial and venous circulatory networks, based on an electrical analogy.
The different compartments are modeled as RLC (resistance, inductance, capacitance) circuits, while cardiac valves are described as diodes.

Among the parameters associated with the $\modEMfom$ model, in this work we focus on the four ones reported in Tab.~\ref{tab:paramsM}.
Similarly, we report in Tab.~\ref{tab:paramsC} the parameters associated with the $\modCirc$ model.

\begin{table}
    \begin{center}
        \begin{tabular}{ rrrl }
            \toprule
            Parameter & Baseline & Unit & Description \\
            \midrule
            $\aXB$         & 160.0 & $\si{\mega\pascal}$ & Cardiomyocytes contractility \\
            $\EPcondf$     & 76.43 & $\si{\milli\meter\per\second}$ & Electrical conductivity along fibers \\
            $\fibersAngle$ & 60.0 & degrees & Fibers angle rotation \\
            $\MstiffPass$  & 0.88 & $\si{\kilo\pascal}$ & Passive stiffness \\
            \bottomrule
        \end{tabular}
        \caption{Parameters $\paramM$ of the $\modEMfom$ model considered in this work and associated baseline values.}
        \label{tab:paramsM}
    \end{center}
\end{table}

\begin{table}
    \begin{center}
        \begin{tabular}{ llrl }
            \toprule
            Parameter & Baseline & Unit & Description \\
            \midrule
            $\Vheart$ & 400.15 & \si{\milli\liter} & Initial blood pool volume of the heart \\
            $\EaLA$, $\EaRA$, $\EaRV$ & 0.07, 0.06, 0.55 & \si{\mmHg \per \milli\liter} & LA/RA/RV active elastance \\
            $\EpLA$, $\EpRA$, $\EpRA$ & 0.18, 0.07, 0.05 & \si{\mmHg \per \milli\liter} & LA/RA/RV passive elastance \\
            $\TCLA$, $\TCRA$, $\TCRV$ & 0.14, 0.14, 0.20 & \si{\second}                 & LA/RA/RV contraction time \\
            $\TRLA$, $\TRRA$, $\TRRV$ & 0.14, 0.14, 0.32 & \si{\second}                 & LA/RA/RV relaxation time \\
            $\VnLA$, $\VnRA$, $\VnRV$ & 4.0, 4.0, 16.0 & \si{\milli\liter}              & LA/RA/RV reference volume \\
            $\tLAV$, $\tRAV$ & 0.16, 0.16 & \si{\second}                                & Left/Right atrioventricular delay \\
            $\Rmin$, $\Rmax$ & 0.0075, 75006.2 & \si{\mmHg \second \per \milli\liter}   & Valve minimum/maximum resistance \\
            $\RarSYS$, $\RvnSYS$ & 0.64, 0.32 & \si{\mmHg \second \per \milli\liter}  & Systemic arterial/venous resistance \\
            $\RarPUL$, $\RvnPUL$ & 0.032, 0.036 & \si{\mmHg \second \per \milli\liter}  & Pulmonary arterial/venous resistance \\
            $\CarSYS$, $\CvnSYS$ & 1.2, 60.0 & \si{\milli\liter \per \mmHg}  & Systemic arterial/venous capacitance \\
            $\CarPUL$, $\CvnPUL$ & 10.0, 16.0 & \si{\milli\liter \per \mmHg}  & Pulmonary arterial/venous capacitance \\
            $\LarSYS$, $\LvnSYS$ & 0.005, 0.0005 & \si{\mmHg \second\squared \per \milli\liter}  & Systemic arterial/venous inductance \\
            $\LarPUL$, $\LvnPUL$ & 0.0005, 0.0005 & \si{\mmHg \second\squared \per \milli\liter}  & Pulmonary arterial/venous inductance \\
            \bottomrule
        \end{tabular}
        \caption{Parameters $\paramC$ of the $\modCirc$ model considered in this work and associated baseline values (LA = left atrium; RA = right atrium; RV = right ventricle).}
        \label{tab:paramsC}
    \end{center}
\end{table}

We report in Tab.~\ref{tab:QoI} the list of all the QoIs, computed from the solution of the $\modEMCfom$ model, and outputs of interest that are used through this paper.
The last column of the table indicates which variables are used for cross-validation during the training phase (see Sec.~\ref{sec:methods:hyperparameters}), those that are used for sensitivity analysis purposes (see Sec.~\ref{sec:methods:GSA}) and those that are used to run Bayesian parameter estimations (see Sec.~\ref{sec:methods:bayesianPE}).

\begin{table}
    \begin{center}
        \begin{tabular}{ llll }
            \toprule
            Parameter & Unit & Description & Usage
            \\
            \midrule
            \multicolumn{4}{l}{\textbf{Left Atrium}}\\
            $\VminLA, \VmaxLA$    & \si{\milli\liter} & End-systolic and end-diastolic volume & GSA \\
            $\pminLA, \pmaxLA$    & \si{\mmHg}        & Minimum and maximum pressure & GSA \\
            \multicolumn{4}{l}{\textbf{Left ventricle}}\\
            $\VLV(t)$       & \si{\milli\liter} & Volume transient & X-validation \\
            $\PLV(t)$       & \si{\mmHg}        & Pressure transient & X-validation \\
            $\VminLV, \VmaxLV$    & \si{\milli\liter} & End-systolic and end-diastolic volume & X-validation, GSA \\
            $\pminLV, \pmaxLV$    & \si{\mmHg}        & Minimum and maximum pressure & X-validation, GSA \\
            $\SVLV$      & \si{\milli\liter} & Stroke volume ($\VmaxLV - \VminLV$) & X-validation, GSA \\
            \multicolumn{4}{l}{\textbf{Right Atrium}}\\
            $\VminRA, \VmaxRA$    & \si{\milli\liter} & End-systolic and end-diastolic volume & GSA \\
            $\pminRA, \pmaxRA$    & \si{\mmHg}        & Minimum and maximum pressure & GSA \\
            \multicolumn{4}{l}{\textbf{Right ventricle}}\\
            $\VminRV, \VmaxRV$    & \si{\milli\liter} & End-systolic and end-diastolic volume & GSA \\
            $\pminRV, \pmaxRV$    & \si{\mmHg}        & Minimum and maximum pressure & GSA \\
            $\SVRV$      & \si{\milli\liter} & Stroke volume ($\VmaxRV - \VminRV$) & GSA \\
            \multicolumn{4}{l}{\textbf{Systemic arterial circulation}}\\
            $\pminARSYS, \pmaxARSYS$ & \si{\mmHg}        & Minimum and maximum pressure & GSA, MCMC \\
            \bottomrule
        \end{tabular}
        \caption{List of QoIs used through this paper, either for cross-validation (X-validation), global sensitivity analysis (GSA) or MCMC based Bayesian parameter estimation (MCMC).}
        \label{tab:QoI}
    \end{center}
\end{table}

To numerically approximate this multiphysics and multiscale model, we adopt the segregated approach proposed in \cite{regazzoni2021em_circulation_pt2}, by which the subproblems are solved sequentially.
For space discretization, we rely on bilinear Finite Elements defined on hexahedral meshes, adopting a different spatial resolution for the electrophysiological and the mechanical variables.
For time discretization, we employ a staggered scheme, where different time steps are used according to the specific subproblem.
Moreover, to avoid the numerical oscillations arising from the mechanical feedback on force generation (that are commonly cured by recurring to a monolithic scheme \cite{levrero2020sensitivity,pathmanathan2010cardiac}), we use the stabilized-staggered scheme that we proposed in \cite{regazzoni2020oscillation}.
This numerical model requires, on a 32 cores computer platform, nearly 4 hours of computational time to simulate a heartbeat.
More details on the numerical discretization and the computational platform employed to generate the training data used in this paper are available in Appendix~\ref{app:params:EM}.

\subsection{Software libraries}

The cardiac electromechanics simulations considered in this paper are performed by means of the \texttt{life\textsuperscript{x}} library\footnote{\url{https://lifex.gitlab.io}}, a high-performance \texttt{C++} platform developed within the iHEART project\footnote{iHEART - An Integrated Heart Model for the simulation of the cardiac function, European Research Council (ERC) grant agreement No 740132, P.I. Prof. A. Quarteroni}.
To train the ANN-based models, we employ the open source MATLAB library \texttt{model-learning}\footnote{\url{https://model-learning.readthedocs.io/}}, that implements the machine learning method proposed in \cite{regazzoni2019modellearning} and used in this paper.
Sensitivity analysis is carried out through the open source Python library \texttt{SALib}\footnote{\url{https://salib.readthedocs.io/}} \cite{Herman2017}.
Finally, for the MCMC based Bayesian parameter estimation we rely on the open source Python library \texttt{UQpy}\footnote{\url{https://uqpyproject.readthedocs.io/}} \cite{UQpy2020}.
To employ the ANN-based models trained with the MATLAB library \texttt{model-learning} within the Python environment of \texttt{SALib} and \texttt{UQpy}, we exploit \texttt{pyModelLearning}, a Python wrapper for the \texttt{model-learning} library, which is available in its GitHub repository\footnote{\url{https://github.com/FrancescoRegazzoni/model-learning}}.

Both the MATLAB codes that are used to train the ANN-based ROMs and the datasets that permit to reproduce the numerical results are publicly available in the online repository accompanying this manuscript\footnote{\url{https://github.com/FrancescoRegazzoni/cardioEM-learning}}.

%% file: parts_results.tex
\section{Results}
\label{sec:results}

\subsection{Trained models}
\label{sec:results:models}

To generate the pressure and volume transients needed to train an ANN-based ROM, we sample the parameter space $\paramSpaceM \times \paramSpaceC$ with a Monte Carlo approach, even if more sophisticated sampling strategies -- such as Latin Hypercube Design -- can be considered as well.
For simplicity, in this stage we only consider a subset of the parameters $\paramC$, selected as the most significant ones, on the basis of a preliminary variance-based global sensitivity analysis, obtained with a version of the closed-loop model in which the LV is also represented by a 0D circuit element (see Appendix~\ref{app:params:ANN}).

We consider two scenarios.
First, we study the variability with respect to a single parameter, namely the active contractility.
Thus, we set $\paramM = [\aXB]$.
Under this setting, we generate 30 numerical simulations through the $\modEMCfom$ model to train a ROM, henceforth denoted by $\modEMredA$.
For this ROM, the remaining parameters (i.e. $\EPcondf$, $\fibersAngle$ and $\MstiffPass$) are kept constant (more precisely, equal to the baseline values of Tab.~\ref{tab:paramsM}).
Then, we consider the full parametric variability (that is, we set $\paramM= [\aXB, \EPcondf, \fibersAngle, \MstiffPass]$), we generate 40 training samples and we train a second ROM, denoted by $\modEMredB$.
All the numerical simulations included in the training set are 5 heartbeats long.
The optimal sets of hyperparameters, tuned through the algorithm of Sec.~\ref{sec:methods:hyperparameters}, are reported in Tab.~\ref{tab:ROMs}.
In both the considered cases, training an ANN-based model takes approximately 18 hours on a single core standard laptop.

\begin{table}
    \begin{center}
        \begin{tabular}{ l|rr|rrrr }
            \toprule
            Trained model & Parameters & Training set size & \multicolumn{4}{c}{Hyperparameters}
            \\
            & $\paramM$ & $\numTrain$ & $\NumANNState$ & $\ANNnumLayers$ & $\ANNnumNeurons$ & $\ANNregularization$
            \\
            \midrule
            $\modEMredA$ & $[\aXB]$                                      & 30 & 2 & 1 & 8  & 0
            \\[1mm]
            $\modEMredB$ & $[\aXB, \EPcondf, \fibersAngle, \MstiffPass]$ & 40 & 1 & 1 & 12 & $0.01$
            \\
            \bottomrule
        \end{tabular}
        \end{center}
        \caption{Optimal sets of hyperparameters for the two trained models $\modEMredA$ and $\modEMredB$.}
        \label{tab:ROMs}
\end{table}

Once trained, the two ROMs ($\modEMredA$ and $\modEMredB$) can be coupled with the circulation model $\modCirc$, thus obtaining the models $\modEMCredA$ and $\modEMCredB$, according to Eq.~\eqref{eqn:model_ANN-0D}.
These two models represent two surrogates of the $\modEMCfom$ model, capable of approximating its output at a dramatically reduced computational cost.
As a matter of fact, numerical simulations with either the $\modEMredA$ or the $\modEMredB$ model take nearly one second of computational time per heartbeat on a single core standard laptop.

To test the accuracy of the $\modEMCredA$ and $\modEMCredB$ models with respect to the $\modEMCfom$ model, we consider a testing dataset, by taking unobserved samples in the parameter space $\paramSpaceM \times \paramSpaceC$.
For both models, we consider 15 testing simulations of the same duration of the ones included in the training set.
In addition, to test the reliability of the ROMs over a longer time horizon than the one considered in the training set, we include in the testing set 5 simulations of double length (i.e. 10 heartbeats).
Then, we compare the simulations obtained with the two ROMs ($\modEMredA$ and $\modEMredB$) with the ones obtained with the $\modEMCfom$ model for the same parameters $\paramM$ and $\paramC$.

The accuracy of the two ROMs is summarized Tabs.~\ref{tab:errors_LV} and \ref{tab:errors_RV}.
Specifically, in Tab.~\ref{tab:errors_LV} we report metrics regarding the ability of the ROMs to correctly predict the function of the LV, that is the chamber surrogated by the ANN-based model.
Specifically, we report the relative errors in $L^2$ norm (i.e. the mean square errors) associated with pressure and volume transients ($\PLV$ and $\VLV$) and the relative errors on some biomarkers of clinical interest (maximum and minimum pressures and volumes).
For the latter, we also report the coefficient of determination $\Rtwo$.
We notice that both model $\modEMCredA$ and model $\modEMCredB$ are able to surrogate model $\modEMCfom$ with remarkably good accuracy.
Model $\modEMCredB$ features slightly larger errors than model $\modEMCredA$; this is not surprising since model $\modEMCredB$ explores a much larger parametric space than model $\modEMCredA$ (four parameters instead of one).
Interestingly, the errors obtained over a long time horizon are very similar to those obtained for simulations of the same duration as those used to train the ANNs.
This demonstrates the reliability of the ROMs in long-term simulations.
Similarly, in Tab.~\ref{tab:errors_RV} we report the errors and $\Rtwo$ coefficients associated with the RV output (for simplicity, we here consider only 5 heartbeats long simulations).
The accuracy obtained in reproducing the RV function is even better than that for the LV, coherently with the fact that the RV is included in the $\modCirc$ model and thus it is not surrogated by the ANN-based model.

\begin{table}
\begin{center}
       \begin{tabular}{ rlrrrrrr }
       \toprule
       & & \multicolumn{6}{c}{\textbf{5 heartbeats}} \\
       & & $\PLV(t)$ & $\VLV(t)$ & $\pminLV$ & $\pmaxLV$ & $\VminLV$ & $\VmaxLV$
       \\
       \multirow{ 2}{*}{$\modEMCredA$ vs $\modEMCfom$} & \multicolumn{1}{l|}{relative error} &
       0.0336 & 0.0090 & 0.0097 & 0.0046 & 0.0139 & 0.0035
       \\
       & \multicolumn{1}{l|}{\Rtwo} &
              &        & 99.691 & 99.864 & 99.896 & 99.948
       \\
       \multirow{ 2}{*}{$\modEMCredB$ vs $\modEMCfom$} & \multicolumn{1}{l|}{relative error} &
       0.0620 & 0.0285 & 0.0517 & 0.0272 & 0.0471 & 0.0127
       \\
       & \multicolumn{1}{l|}{\Rtwo} &
              &        & 94.370 & 95.302 & 95.942 & 97.061
       \\
       \midrule
       & & \multicolumn{6}{c}{\textbf{10 heartbeats}} \\
       & & $\PLV(t)$ & $\VLV(t)$ & $\pminLV$ & $\pmaxLV$ & $\VminLV$ & $\VmaxLV$
       \\
       \multirow{ 2}{*}{$\modEMCredA$ vs $\modEMCfom$} & \multicolumn{1}{l|}{relative error} &
       0.0293 & 0.0071 & 0.0113 & 0.0037 & 0.0096 & 0.0031
       \\
       & \multicolumn{1}{l|}{\Rtwo} &
              &        & 99.924 & 99.980 & 99.851 & 99.944
       \\
       \multirow{ 2}{*}{$\modEMCredB$ vs $\modEMCfom$} & \multicolumn{1}{l|}{relative error} &
       0.0631 & 0.0265 & 0.0442 & 0.0147 & 0.0382 & 0.0122
       \\
       & \multicolumn{1}{l|}{\Rtwo} &
              &        & 92.227 & 99.957 & 99.229 & 99.063
       \\
       \bottomrule
       \end{tabular}
       \caption{Testing errors and \Rtwo coefficients on the LV outputs obtained with the two models $\modEMredA$ and $\modEMredB$, for 5 heartbeats long (top) and 10 heartbeats long (bottom) simulations.}
       \label{tab:errors_LV}
\end{center}
\end{table}

\begin{table}
\begin{center}
       \begin{tabular}{ rlrrrrrr }
       \toprule
       & & \multicolumn{6}{c}{\textbf{5 heartbeats}} \\
       & & $\PRV(t)$ & $\VRV(t)$ & $\pminRV$ & $\pmaxRV$ & $\VminRV$ & $\VmaxRV$
       \\
       \multirow{ 2}{*}{$\modEMCredA$ vs $\modEMCfom$} & \multicolumn{1}{l|}{relative error} &
       0.0048 & 0.0035 & 0.0015 & 0.0027 & 0.0028 & 0.0004
       \\
       & \multicolumn{1}{l|}{\Rtwo [\%]} &
              &        & 100.000 & 99.993 & 99.998 & 100.000
       \\
       \multirow{ 2}{*}{$\modEMCredB$ vs $\modEMCfom$} & \multicolumn{1}{l|}{relative error} &
       0.0069 & 0.0040 & 0.0029 & 0.0079 & 0.0072 & 0.0069
       \\
       & \multicolumn{1}{l|}{\Rtwo [\%]} &
              &        & 99.994 & 99.807 & 99.819 & 99.997
       \\
       \bottomrule
       \end{tabular}
       \caption{Testing errors and \Rtwo coefficients on the RV outputs obtained with the two models $\modEMredA$ and $\modEMredB$, for 5 heartbeats long simulations.}
       \label{tab:errors_RV}
\end{center}
\end{table}

\begin{figure}
       \centering
       \includegraphics{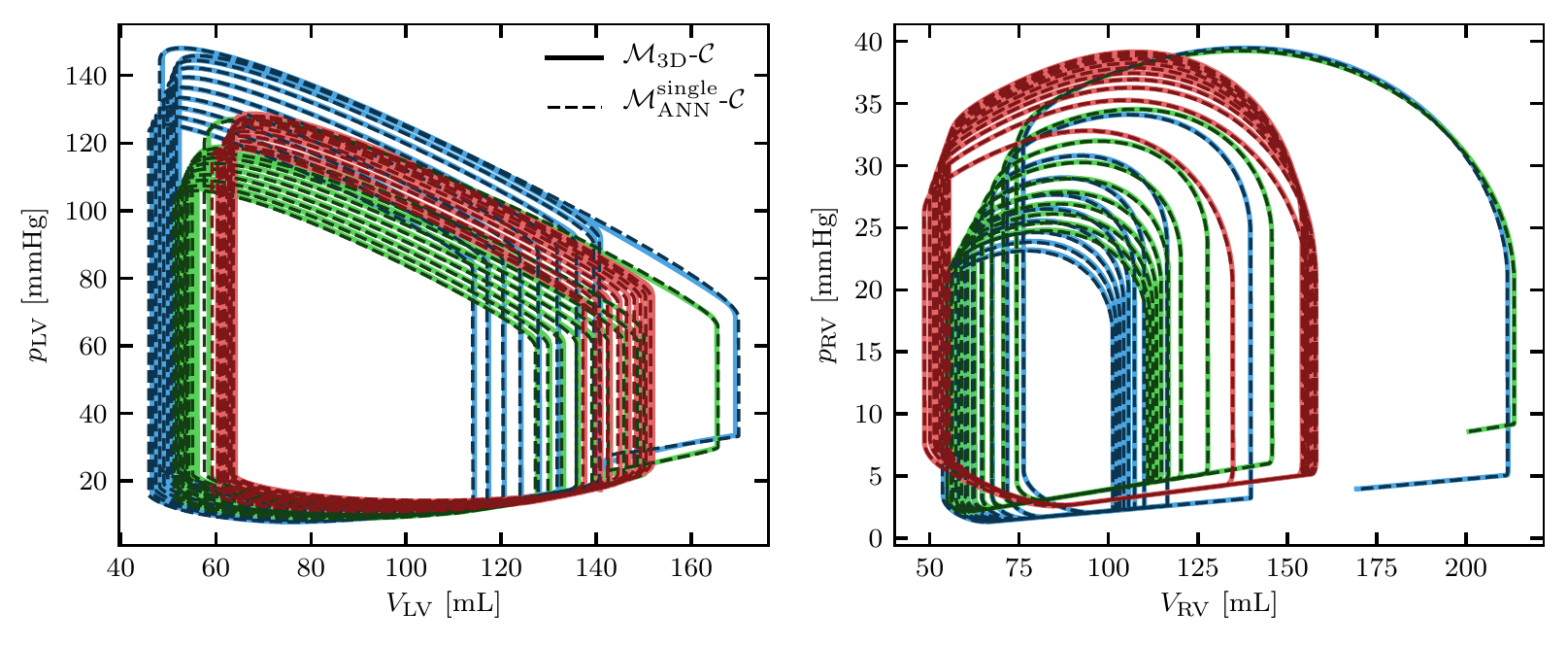}
       \includegraphics{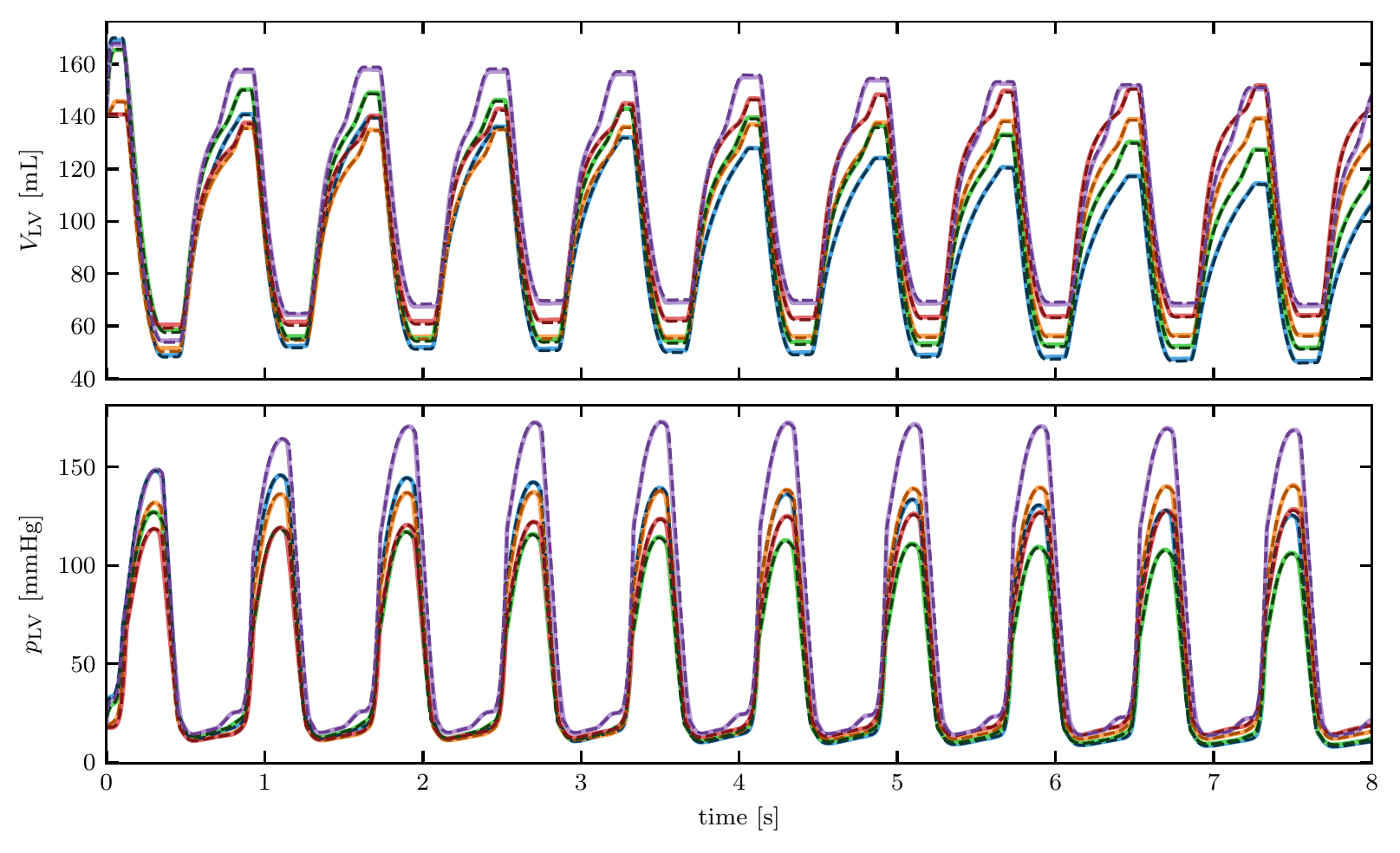}
       \caption{
       Pressure and volume transients obtained with the $\modEMCredA$ (dashed lines), compared to those obtained with the $\modEMCfom$ model (solid lines).
       The different colors correspond to different samples of the testing set.
       For the sake of clarity, only three samples are shown in the first row.
       }
       \label{fig:transients_a_XB}
\end{figure}

\begin{figure}
       \centering
       \includegraphics{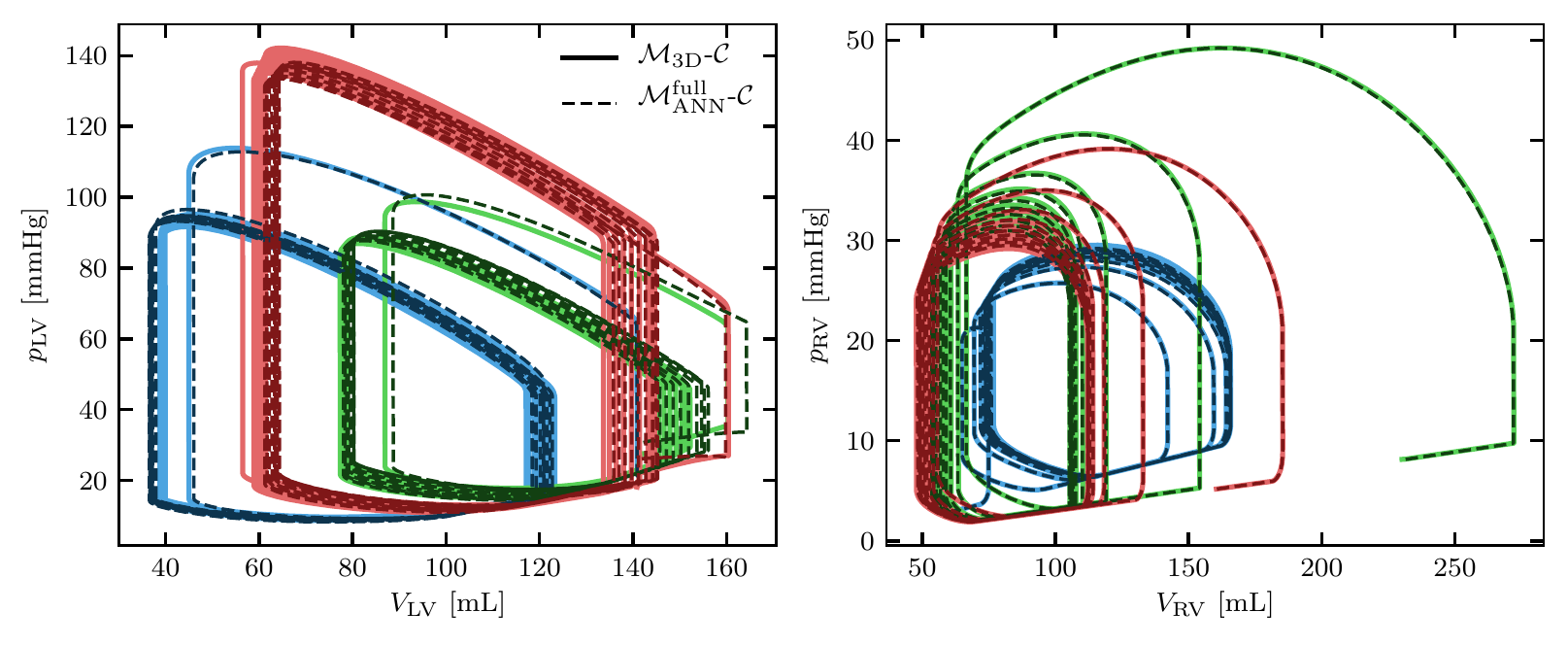}
       \includegraphics{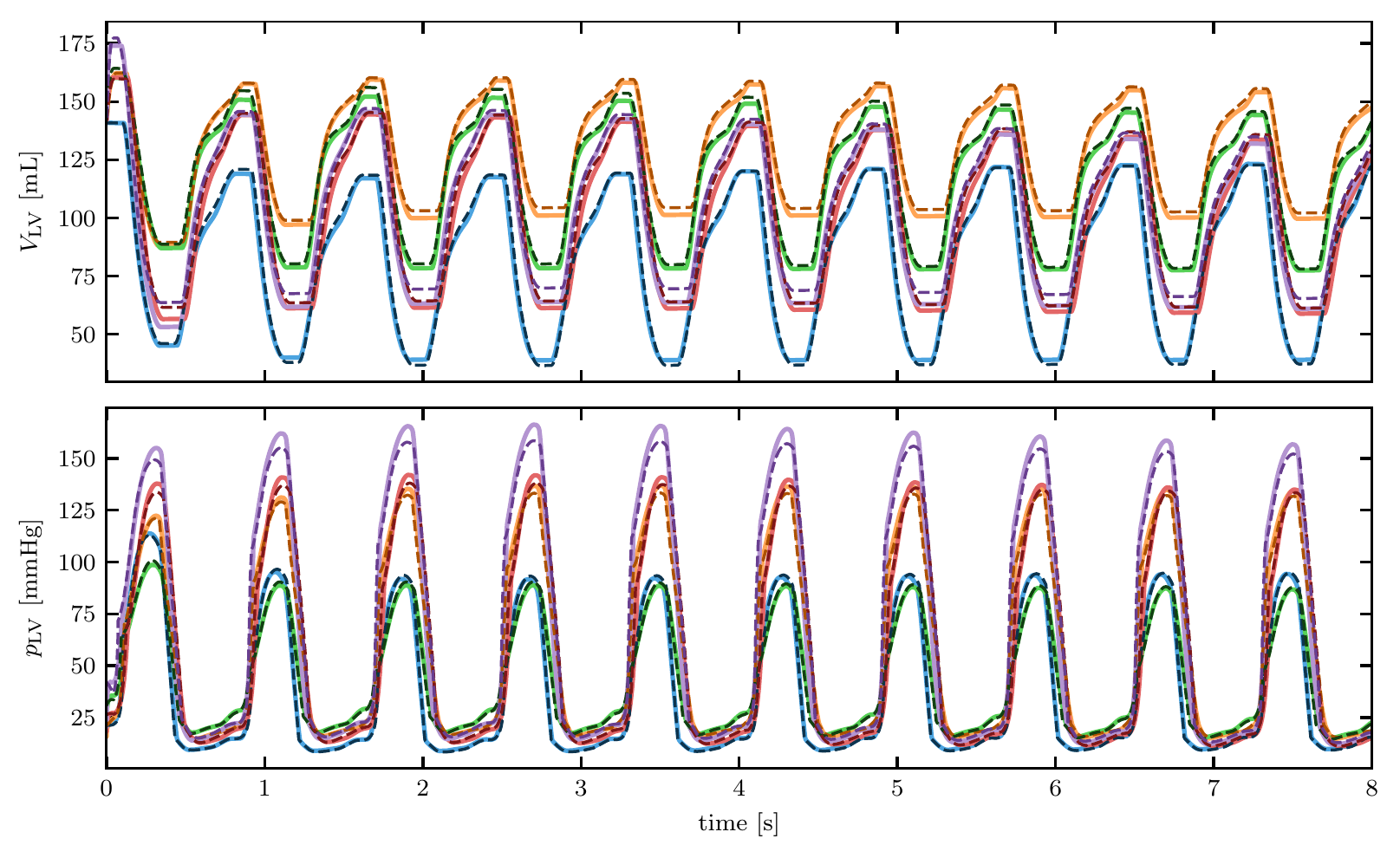}
       \caption{
       Pressure and volume transients obtained with the $\modEMCredB$ (dashed lines), compared to those obtained with the $\modEMCfom$ model (solid lines).
       The different colors correspond to different samples of the testing set.
       For the sake of clarity, only three samples are shown in the first row.
       }
       \label{fig:transients_full}
\end{figure}

\begin{figure}
       \centering
       {\large Model $\modEMCfom$ vs $\modEMCredA$}
       \\
       \includegraphics{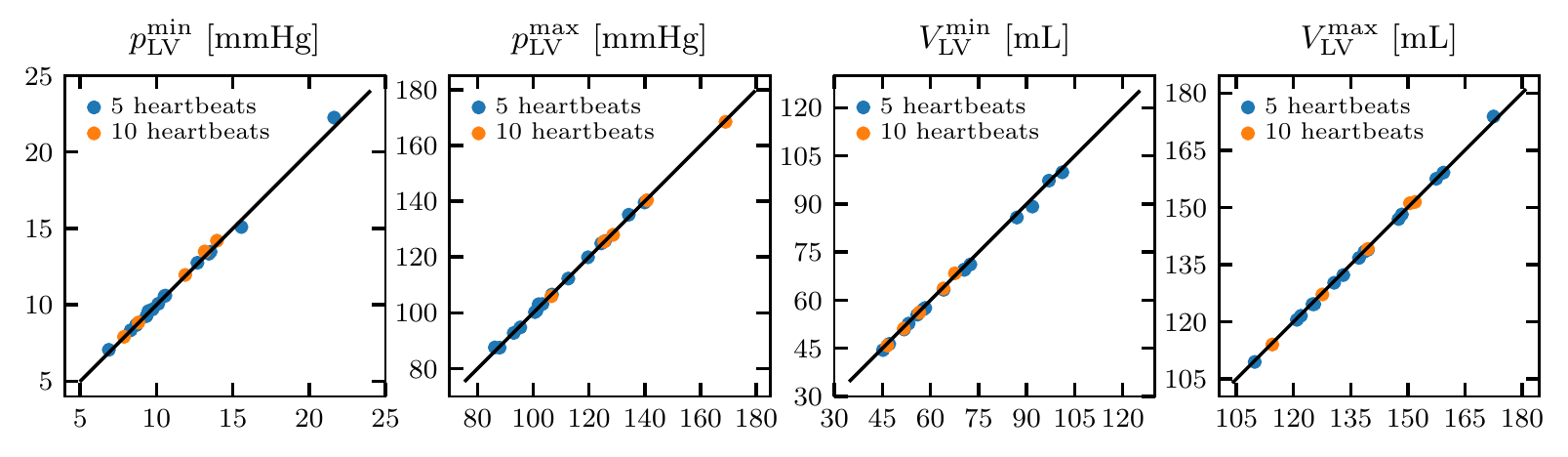}
       \\[3mm]
       {\large Model $\modEMCfom$ vs $\modEMCredB$}
       \\
       \includegraphics{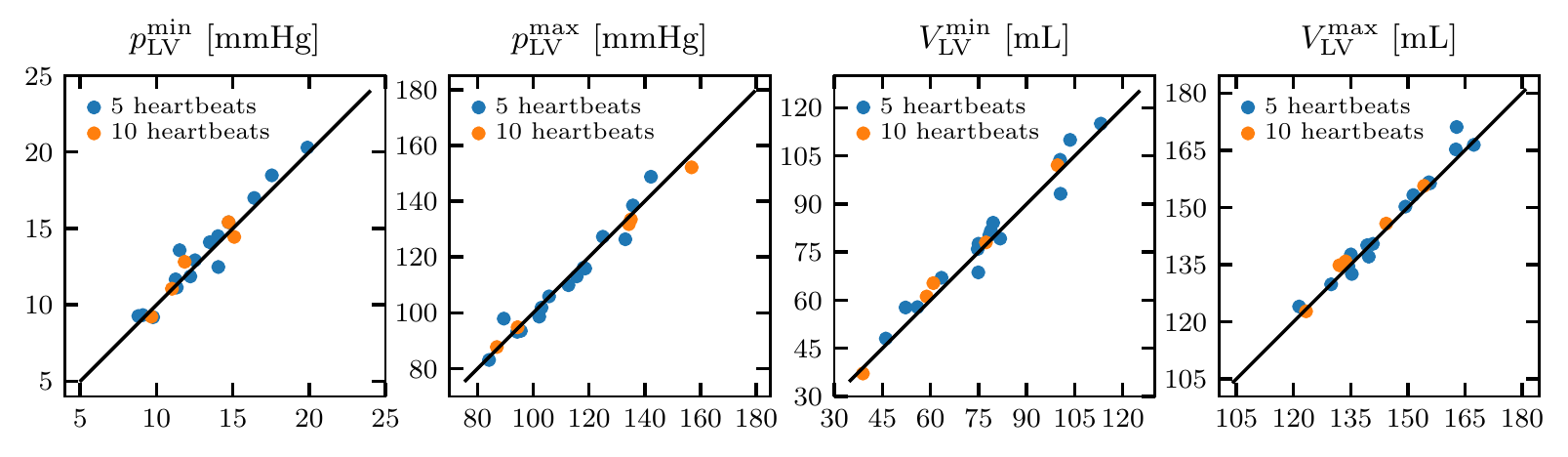}
       \caption{
       LV biomarkers obtained with the $\modEMCredA$ and $\modEMCredB$ models versus those obtained with the $\modEMCfom$ model in the testing set.
       The different marker colors are associated with 5 heartbeats and 10 heartbeats long simulations, respectively.
       }
       \label{fig:biomarkers}
\end{figure}

In Figs.~\ref{fig:transients_a_XB} and \ref{fig:transients_full} we compare 10 heartbeats long transients obtained with the $\modEMCredA$ and $\modEMCredB$ models, respectively, to those obtained with the $\modEMCfom$ model.
All these transients were not used to train the ANN-based models (that is, they belong to the testing set).
In Fig.~\ref{fig:biomarkers} we show the LV biomarkers predicted by the $\modEMCredA$ and $\modEMCredB$ models versus those predicted by the $\modEMCfom$ model.

\subsection{Coupling the electromechanical reduced-order model with different circulation models}
\label{sec:results:windkessel}

To generate the data used to train the $\modEMred$ model, we employ the $\modEMCfom$ coupled model (that is, we couple the $\modEMfom$ model with the circulation model $\modCirc$).
However, the ANN-based ROM $\modEMred$ surrogates the $\modEMfom$ model regardless of its coupling with a specific circulation model.
In fact, thanks to Eq.~\eqref{eqn:model_ANN-0D}, the trained $\modEMred$ model can be coupled with circulation models that are different from the one used to generate the training data.
We remind that predictions are reliable only if pressure and volume values are inside the ranges explored during the training phase.

We demonstrate the flexibility of our approach by coupling the $\modEMred$ model with a pressure-volume closure relationship that is different from the closed-loop circulation model $\modCirc$ used during the training phase (see Sect.~\ref{sec:methods:EM_model}).
Specifically, we consider a circulation model $\modCircW$ with a windkessel type relationship during the ejection phase and a linear pressure ramp during the filling phase \cite{regazzoni2019mor-sarcomeres}.
In spite of the different nature of the two circulation models (the one used during the training and the one used during the testing), the ANN-based ROM trained with the $\modCirc$ model proves to be reliable also when it is coupled with the $\modCircW$ model.
Indeed, as shown in Tab.~\ref{tab:errors_windkessel}, the results obtained by the $\modEMCWred$ coupled model approximate those of the $\modEMCWfom$ coupled model with an excellent accuracy.
These errors are indeed comparable to the ones obtained by surrogating the $\modEMCred$ model with the $\modEMCfom$ model (see Tab.~\ref{tab:errors_LV}).

\begin{table}
       \begin{center}
              \begin{tabular}{ rlrrrrrr }
              \toprule
              & & $\PLV(t)$ & $\VLV(t)$ & $\pminLV$ & $\pmaxLV$ & $\VminLV$ & $\VmaxLV$
              \\
              $\modEMCWredA$ vs $\modEMCWfom$ & \multicolumn{1}{l|}{relative error} &
              0.0459 & 0.0518 & 0.0065 & 0.0004 & 0.0080 & 0.0009
              \\
              $\modEMCWredB$ vs $\modEMCWfom$ & \multicolumn{1}{l|}{relative error} &
              0.0510 & 0.0219 & 0.0116 & 0.0027 & 0.0110 & 0.0063
              \\
              \bottomrule
              \end{tabular}
              \caption{Testing errors on the LV outputs obtained with the two models $\modEMredA$ and $\modEMredB$ coupled with the $\modCircW$ model for 5 heartbeats long simulations.}
              \label{tab:errors_windkessel}
       \end{center}
\end{table}

\subsection{Global sensitivity analysis}
\label{sec:results:SA}

Once we have verified that models $\modEMCredA$ and $\modEMCredB$ are able to reproduce the outputs of model $\modEMCfom$ with high accuracy, we use them to perform a global senstivity analysis.
The aim is to determine which of the parameters of the circulation model ($\paramC$) and the electromechanical model ($\paramM$) contribute the most to the determination of a list of outputs of interest, the so-called QoIs (see Tab.~\ref{tab:QoI}).
For this purpose, we compute the Sobol indices $\SobolFirst{i}{j}$ and $\SobolTotal{i}{j}$, as described in Sec.~\ref{sec:methods:GSA}.
When replacing model $\modEMCfom$ with a ROM, we can only estimate sensitivity indices with respect to the parameters $\paramM$ that were considered during the training phase.
On the other hand, we can compute sensitivity indices with respect to all parameters $\paramC$ of the circulation model, even those that were not varied during training.
It would even be possible -- at least in principle -- to consider a circulation model different from the one used to generate the training data.
In fact, as pointed out in Sec.~\ref{sec:methods:ROM}, the model $\modEMredA$ and $\modEMredB$ surrogate the models $\modEMfom$ independently of the circulation model to which it is coupled with.

\begin{figure}
       \centering
       \includegraphics{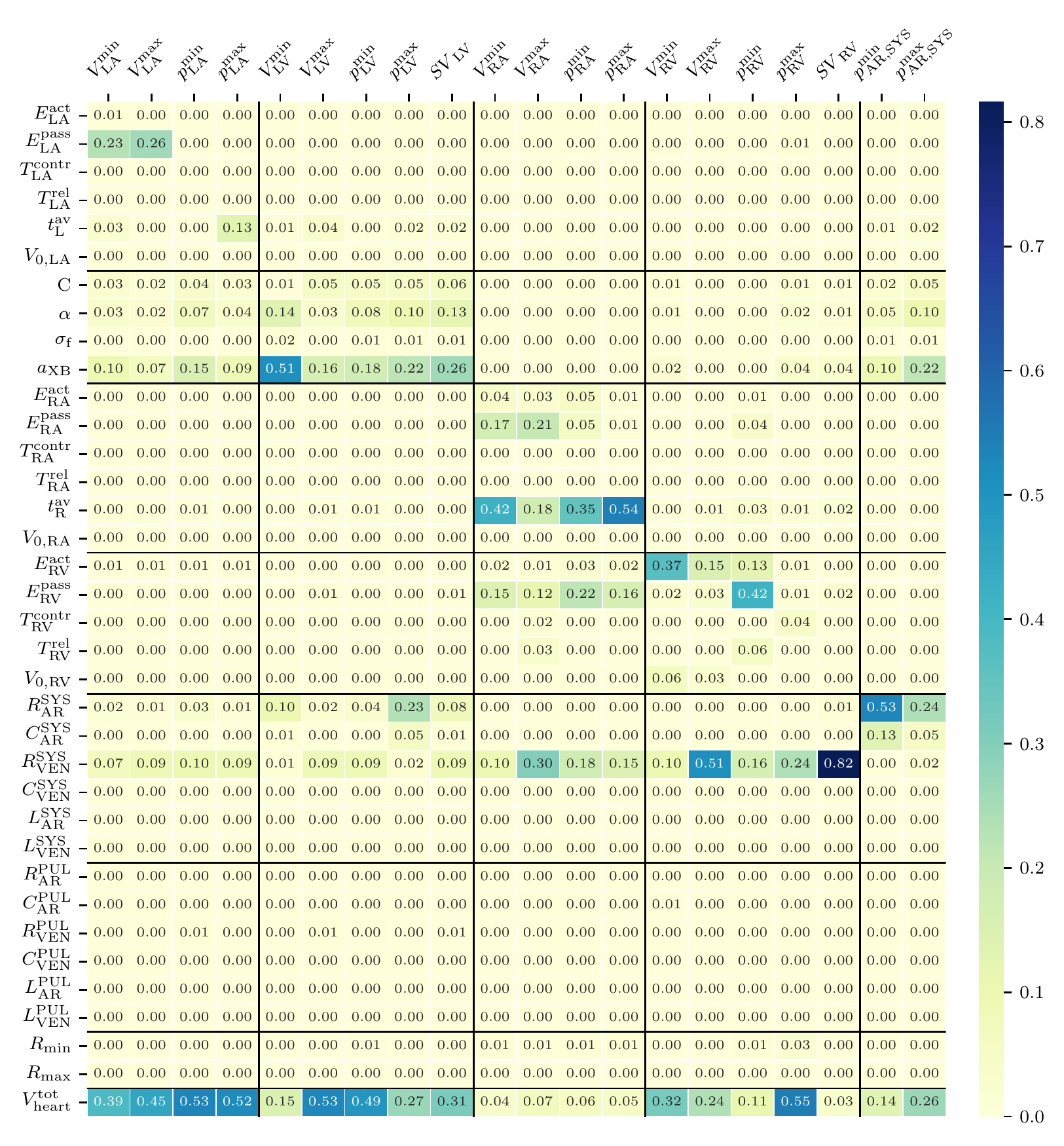}
       \caption{
       First-order Sobol indices $\SobolFirst{i}{j}$ computed by exploiting the $\modEMCredB$ model.
       Each row corresponds to a parameter of either the electromechanical model (i.e. $\paramM$, see Tab.~\ref{tab:paramsM}) or the circulation model (i.e. $\paramC$, see Tab.~\ref{tab:paramsC}).
       Each column corresponds to a QoI (i.e. $\qoi$, see Tab.~\ref{tab:QoI}).
       Both parameters and QoIs are split into a number of groups, separated by a black solid line.
       Specifically, from left to right, we list QoIs referred to LA, LV, RA, RV and systemic circulation.
       Similarly, from top to bottom, we list parameters associated with LA, LV, RA, RV, systemic circulation, pulmonary circulation, valves and blood total volume.
       }
       \label{fig:SA_all_S1}
\end{figure}

\begin{figure}
       \centering
       \includegraphics{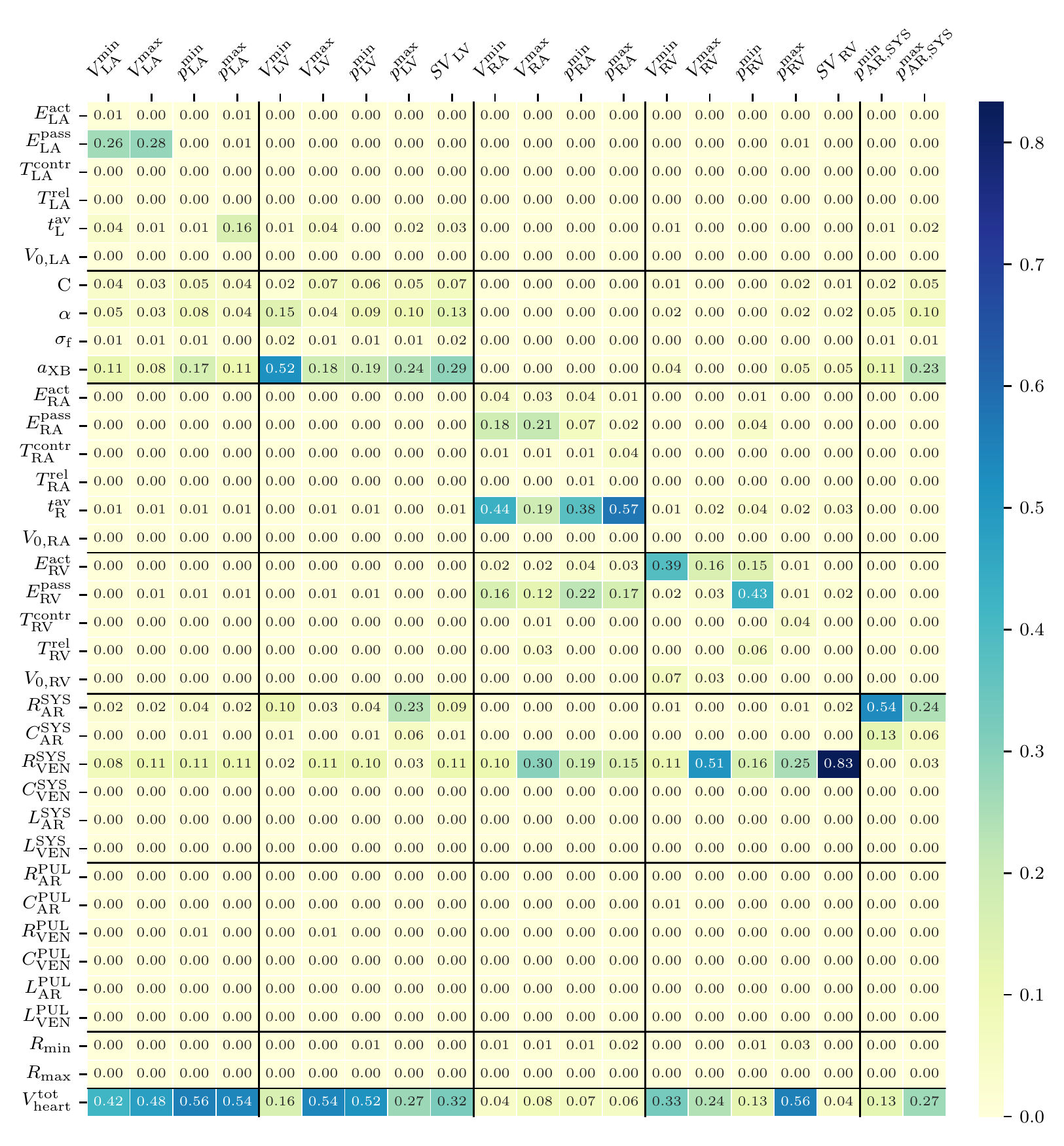}
       \caption{
       Total-effect Sobol indices $\SobolTotal{i}{j}$ computed by exploiting the $\modEMCredB$ model.
       For a description of the figure see caption of Fig.~\ref{fig:SA_all_S1}.
       }
       \label{fig:SA_all_ST}
\end{figure}

We only report the results obtained by means of the most complete ROM, that is $\modEMCredB$.
The results are shown in Figs.~\ref{fig:SA_all_S1} and \ref{fig:SA_all_ST}, respectively.
First, we notice that the $\SobolFirst{i}{j}$ and $\SobolTotal{i}{j}$ indices have only small differences from each other.
This means that the interaction among the different parameters is less significant, in the determination of the QoIs, than the variation of the individual parameters.
Furthermore, we note that, as expected, the QoIs associated with a given chamber or compartment are mostly determined by the parameters associated with the same region of the cardiovascular network.
However, there are some important exceptions.
Indeed, the venous resistance of the systemic circulation $\RvnSYS$ has a strong impact on the right heart (RA and RV), i.e. the part of the network located immediately upstream.
In addition, the systemic arterial resistance $\RarSYS$ significantly impacts the maximum LV pressure.
In fact, this parameter is known to contribute in the determination of the so-called afterload \cite{tortora2008}.
We also note that total circulating blood volume $\Vheart$ has a major impact on almost all biomarkers.
Conversely, the parameters associated with the pulmonary circulation network have a minimal impact.
In addition, parameters describing the resistance of opened and closed valves also have a very little impact.
This is an interesting result, since it shows that these parameters, chosen as a very low and very high value respectively (since for reasons of numerical stability they cannot be set equal to zero and infinity), have virtually no impact (at least within the ranges considered here) on the output quantities of biomechanical interest.

Regarding individual compartments, we note that variability can be explained by only a few parameters.
Specifically, atria are mainly influenced by passive stiffness and atrioventricular delay, whereas for the RV the most relevant parameters are active and passive stiffness.
As expected, the biomarkers associated with the LV -- the only chamber included in the $\modEMfom$ model -- are mainly influenced by the parameters $\paramM$.
Among these, the most influential one is the active contractility $\aXB$, followed the fibers orientation $\fibersAngle$ and, to a lesser extent, the passive stiffness $\MstiffPass$.
Finally, the electrical conductivity $\EPcondf$ has a minimal impact on the biomarkers under consideration.

It is important to note that Sobol indices are affected by the amplitude of the ranges in which the parameters are varied.
In particular, the wider the range associated with a parameter, the greater the associated Sobol indices will be, as the parameter in question potentially generates greater variability in the QoI.
Therefore, the results shown here are valid for the specific ranges we used, which are reported in Appendix~\ref{app:params:GSA}.

\subsection{Bayesian parameter estimation}
\label{sec:results:PE}

\MS{In this section we present a further practical use of the cardiac electromechanics ROM presented in this paper.}
In particular, we show that the ROM can be used to enable Bayesian parameter estimation for the $\modEMCfom$ model, which, due to the prohibitive computational cost, would not be affordable without the use of a ROM.

First, we consider a prescribed value for the parameter vector $\param = (\paramM, \paramC)$ (specifically, we employ the values reported in Tabs.~\ref{tab:paramsM} and \ref{tab:paramsC}) and we run a simulation using the $\modEMCfom$ model.
Based on the output of this simulation, we consider a couple of QoIs, consisting of minimum and maximum arterial pressure (i.e. we set $\qoi = (\pminARSYS, \pmaxARSYS)$).
The choice of these two QoIs is motivated by the fact that they are two variables that can be measured non invasively, and in fact they are often monitored in clinical routine.
We reconstruct the value of a couple of parameters, namely the active contractility $\aXB$ and the systemic arterial resistance $\RarSYS$, assuming known and keeping fixed the values of the remaining parameters.
\MS{Indeed, we aim at demonstrating that the ROM we propose is suitable for estimating parameters of both the $\modEMfom$ model (such as $\aXB$) and of the $\modCirc$ model (such as $\RarSYS$).}
Specifically, in this section we rely on the $\modEMCredA$ model.

To mimic the presence of measurement errors, we artificially add to the exact values of the QoIs $\pminARSYS$ and $\pmaxARSYS$ a synthetic noise, with increasing magnitude.
Specifically, we add an artificial noise by sampling from independent Gaussian variables with zero mean and with variance $\NoiseMagnEXP^2$.
We consider three cases: $\NoiseMagnEXP^2 = 0$ (i.e. no noise), $\NoiseMagnEXP^2 = \SI{0.1}{\mmHg\squared}$ and $\NoiseMagnEXP^2 = \SI{1}{\mmHg\squared}$.

As described in Sec.~\ref{sec:methods:bayesianPE}, we use the MCMC method to derive a Bayesian estimate of the parameters value from (possibly noisy) measurements of the minimum and maximum arterial pressure.
For both unknown parameters we employ a non-informative prior, that is a uniform distribution on the ranges used to train the ROM ($\aXB \in [80, 320] \, \si{\mega\pascal}$ and $\RarSYS \in [0.4, 1.2] \, \si{\mmHg \second \per \milli\liter}$).
According to Sec.~\ref{sec:methods:bayesianPE}, we set $\NoiseCov = \NoiseCovROM + \NoiseCovEXP$, where the experimental measurement error covariance is given by $\NoiseCovEXP = \NoiseMagnEXP^2 \, \mathbb{I}_2$ ($\mathbb{I}_2$ being the 2-by-2 identity matrix) and where the ROM approximation error covariance is estimated from its statistical distribution on the validation set as $\NoiseCovROM = \SI{0.2}{\mmHg\squared} \, \mathbb{I}_2$.
More details on the MCMC setup are available in Appendix~\ref{app:params:bayesianPE}.

In Fig.~\ref{fig:MCMC} we show the posterior distribution $\piPost$ on the parameters pair ($\aXB$, $\RarSYS$) obtained for the three noise levels considered.
We identify with a red line the 90\% credibility region, that is the region in the parameter space with largest posterior probability such that it covers 90\% of $\piPost$.
We notice that for each value of noise, the credibility region contains the exact value of the parameters (namely $\aXB = \SI{160}{\mega\pascal}$ and $\RarSYS = \SI{0.64}{\mmHg \second \per \milli\liter}$), represented by a red star.
As expected, for larger values of noise, the size of the credibility region increases (that is, the estimate is more uncertain).
As a matter of fact, an advantage of Bayesian parameter estimation methods, compared to deterministic methods, is their ability of quantifying the uncertainty associated with the parameters estimate.
A further feature of Bayesian methods stands in capturing correlations among the estimated parameters.
\MS{Indeed, this aspect emerges clearly because of the oblique shape of the credibility regions.}
This is due to the fact that an increase in $\aXB$ or a decrease in $\RarSYS$ leads to similar changes in terms of the measured QoIs ($\pminARSYS$ and $\pmaxARSYS$), thus making their posterior distributions highly correlated.

\begin{figure}
       \centering
       \includegraphics{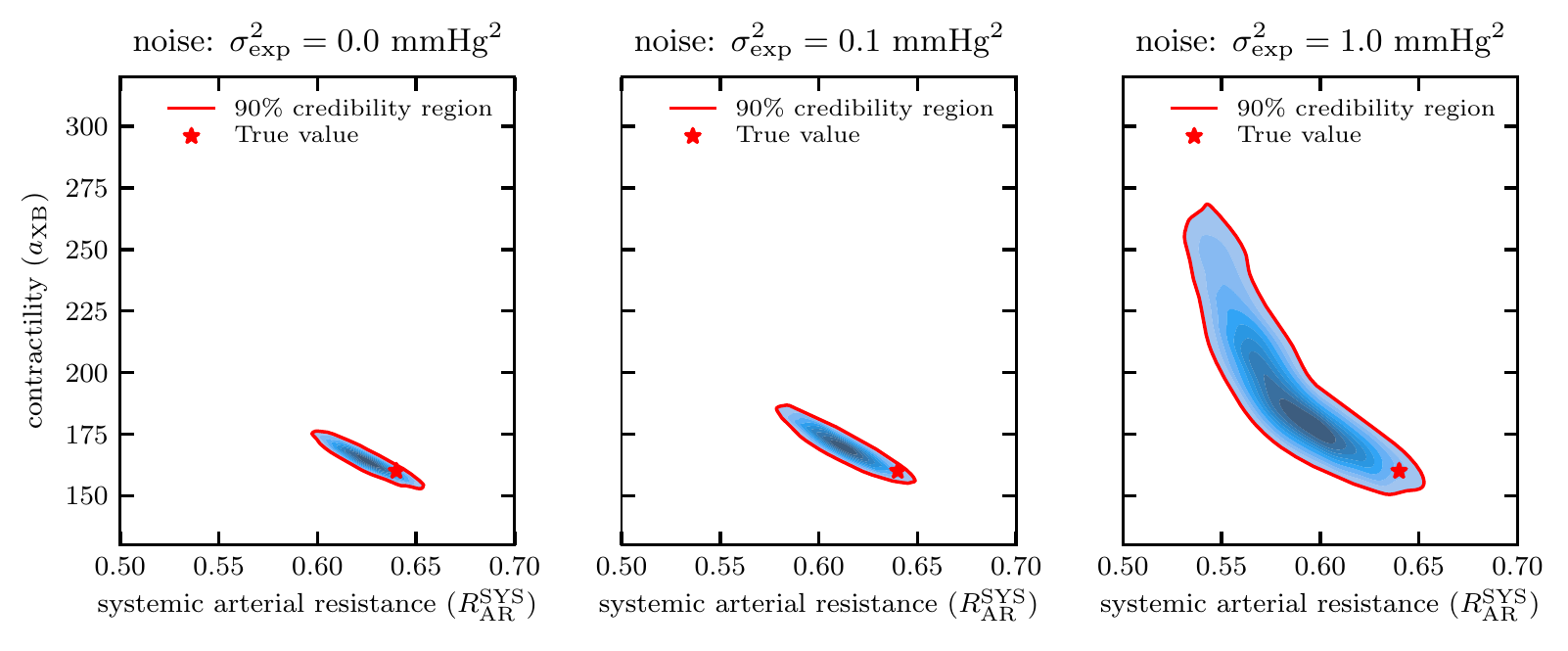}
       \caption{
       \MS{Output of the Bayesian estimation presented in Sec.~\ref{sec:results:PE}.}
       We depict the posterior distribution $\piPost$, estimated by means of the MCMC method, for $\NoiseMagnEXP^2 = 0$ (left), $\NoiseMagnEXP^2 = \SI{0.1}{\mmHg\squared}$ (middle) and $\NoiseMagnEXP^2 = \SI{1}{\mmHg\squared}$ (right).
       The red lines show the 90\% credibility regions, while the red stars represent the exact value of the unknown parameters $\aXB$ and $\RarSYS$.
       }
       \label{fig:MCMC}
\end{figure}

%% file: parts_discussion.tex
\section{Discussion}
\label{sec:discussion}

\begin{figure}
	\centering
	\includegraphics[width=\textwidth]{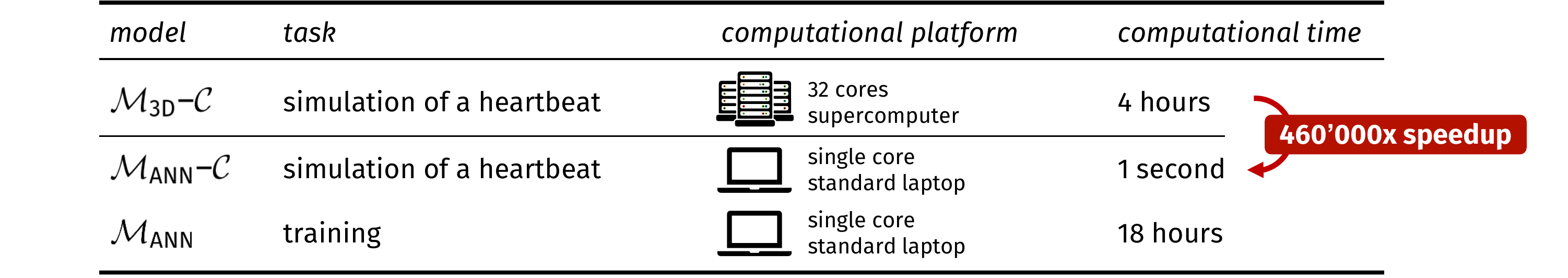}
	\caption{Summary of the computational times associated with FOM and ROM numerical simulations and with training of the ANN-based ROM.}
	\label{fig:computational_times_base}
\end{figure}

In this paper, we develop a machine learning method to build ANN-based ROMs of 3D cardiac \MS{electromechanical} models.
Thanks to the reduced-order differential equations learned through our algorithm, it is possible to approximate the cardiac dynamics in terms of pressure and volume transients with great fidelity and with a huge computational saving.
As a matter of fact, once trained, the ANN-based ROM permits to simulate a heartbeat virtually in real time (about one second of numerical simulation for a heartbeat on a standard laptop), when the original electromechanical model requires about four hours per heartbeat on a supercomputer with 32 cores.
By taking into account the number of cores, our ANN-based ROM brings the impressive speedup of 460'000x (see Fig.~\ref{fig:computational_times_base}).
Clearly, a fair evaluation of the computational saving should also consider the time required to generate the training dataset and to train the model.
To this aim, we consider two test cases, corresponding to the examples of global sensitivity analysis and Bayesian parameter estimation presented in Secs.~\ref{sec:methods:GSA} and \ref{sec:methods:bayesianPE}, respectively.
In both cases, we consider 5 computer nodes with 32 cores each available, and we compare the total computational times using either the $\modEMCfom$ or the $\modEMCred$ model with these computational resources.

\begin{figure}
	\centering
	\includegraphics[width=\textwidth]{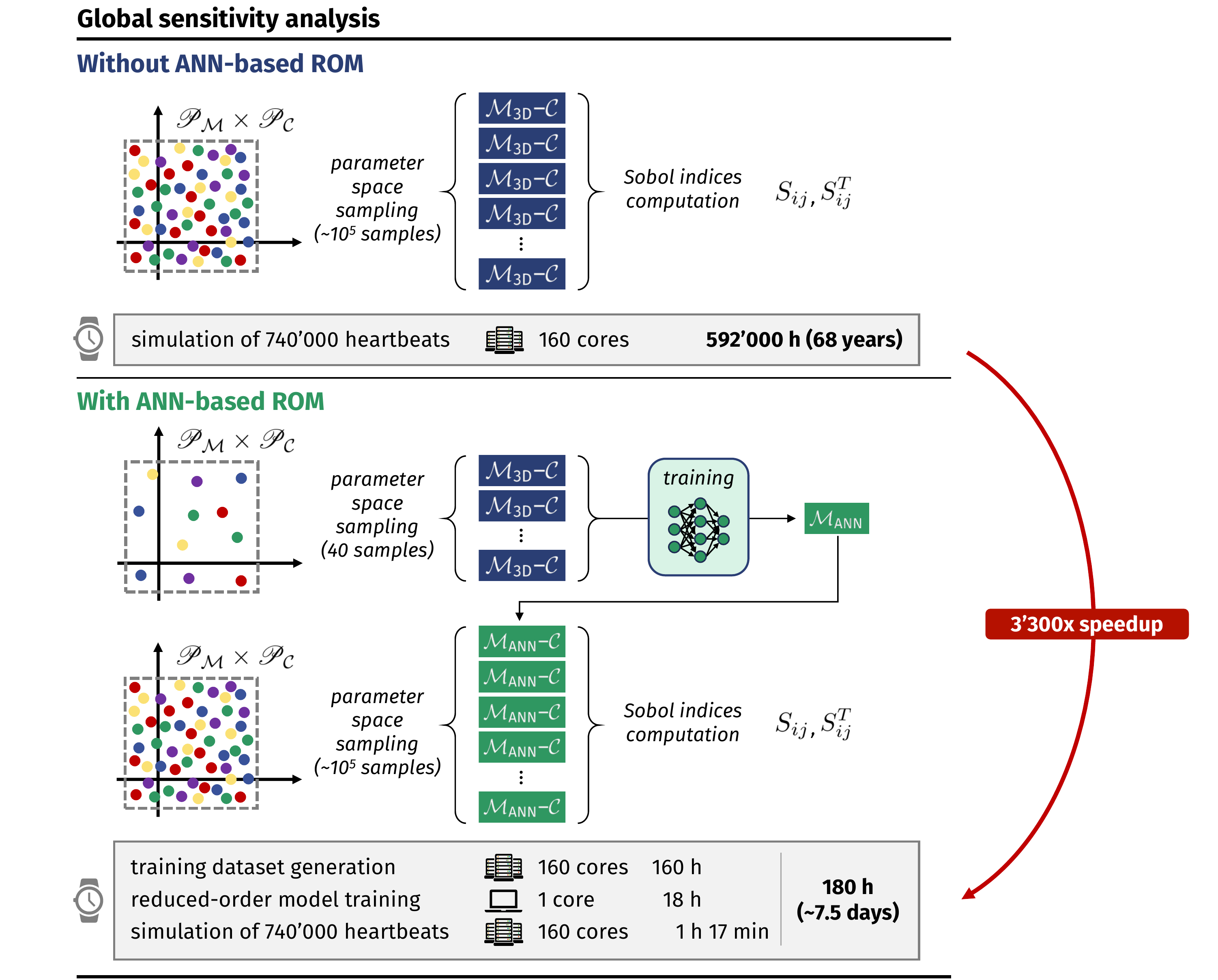}
	\caption{
    Summary of the computational times associated with global sensitivity analysis.
	}
	\label{fig:times_flow_SA}
\end{figure}

\begin{figure}
	\centering
	\includegraphics[width=\textwidth]{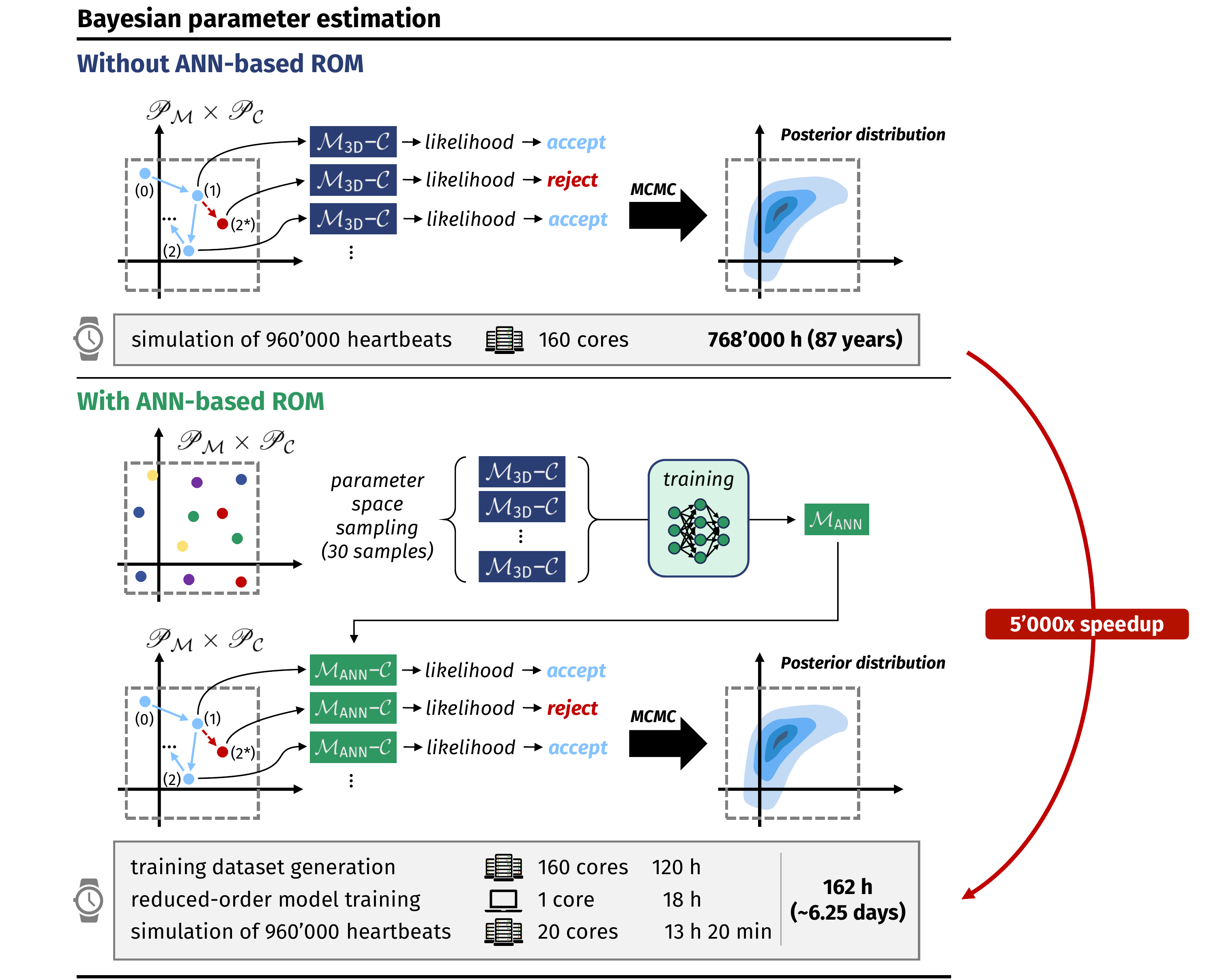}
	\caption{
    Summary of the computational times associated with Bayesian parameter estimation.
	}
	\label{fig:times_flow_BPE}
\end{figure}

The global sensitivity analysis presented in Sec.~\ref{sec:methods:GSA} requires the numerical simulation of 74'000 parameter cases.
On average, the limit cycle is reached in 10 heartbeats, for a total of 740'000 heartbeats to be simulated.
Running this huge amount of numerical simulations with the FOM $\modEMCfom$ model would not be feasible, as it would require about 68 years of uninterrupted use of the 5 nodes endowed with 32 cores.
On the contrary, by virtue of our machine learning method, we are able to perform a global sensitivity analysis, albeit with a small approximation (see Tabs.~\ref{tab:errors_LV} and \ref{tab:errors_RV}) in the numerical results, in 7.5 days (6.7 days to generate the training dataset, 18 hours to train the model and 1 hour and 17 minutes to perform the sensitivity analysis using the $\modEMCred$ model).
Therefore, taking into account the time required to build the ROM, our approach yields a 3'300x speedup (see Fig.~\ref{fig:times_flow_SA}).
Bayesian parameter estimation, on the other hand, requires the numerical simulation of approximately 960'000 heartbeats.
In this case, despite using only 20 cores to run the MCMC, we are able to obtain a result in just 6 days and 8 hours (including generation of the training dataset), when with the FOM $\modEMCfom$ model it would take more than 87 years.
We thus obtain an overall speedup of 5'000x (see Fig.~\ref{fig:times_flow_BPE}).
Moreover, we notice that in case we need to execute a Bayesian calibration for different data, it is not necessary to repeat the training of the ANN-based model, but it is sufficient to re-run only the MCMC algorithm, which takes -- thanks to our ROM -- only 13 hours and 20 minutes.
Finally, we remark that the cost required to train the ANN or to perform numerical simulations with the $\modEMCred$ model does not depend on the specific $\modEMfom$ model at hand, as it is only based on the generated pressure and volume transients.
In particular, it is not expected to raise if the biophysical detail or the number of degrees of freedom of the computational mesh increase.

The ROM proposed in this paper allows, as shown in Sec.~\ref{sec:results}, to tackle many-query problems, which would otherwise be unaffordable due to their computational cost.
Indeed, sensitivity analysis on cardiac models have mainly addressed, at present, single cell models \cite{sher2013local,tondel2015quantifying,pathmanathan2019comprehensive}.
Sensitivity studies on 3D models of cardiac electromechanics have been mainly limited to simplified geometries \cite{Hurtado2017}, or to a few parameters, varying one parameter at a time (that is, \textit{local} sensitivity analysis) \cite{land2012analysis,Marx2020}.
To the best of our knowledge, the only global sensitivity study on a 3D electromechanical model was only recently performed in \cite{levrero2020sensitivity}; however, the huge computational cost associated with the numerical approximation of the forward model dictated the use of an idealized ellipsoidal geometry and of single heartbeat simulations, thus without letting the system reaching a periodic limit cycle.
This leads to results with limited accuracy.
Furthermore, parameter estimation under uncertainty using 3D models of cardiac electromechanics is hard to attain, due to the computational cost that would be required by Bayesian methods such as MCMC.
Indeed, parameters are typically estimated sequentially by using minimization methods \cite{land2012analysis,wang2009modelling,finsberg2018estimating} or ad hoc developed heuristics \cite{Marx2020, costa2013automatic}.
An alternative family of methods, able to provide an estimate of uncertainty, is that of sequential methods of the Kalman filter type (see e.g. \cite{xi2011myocardial,Marchesseau2013}), which, however, would require an invasive implementation in the electromechanical model.

Recently, a number of surrogate models of cardiac electromechanics have been built by using machine learning algorithms and they are often called emulators \cite{di2018gaussian,dabiri2019prediction,Longobardi2020,Cai2021}.
These emulators are based on a collection of pre-computed numerical simulations obtained by sampling the parameter space, similarly to what has been done in this work.
However, the approach behind these emulators is very different from the one we follow.
These emulators are in fact maps that fit the parameters-to-QoIs map ($\forward \colon \param \mapsto \qoi$) in a static manner.
On the contrary, with our approach, the ROM makes it possible to perform a real numerical simulation of the cardiac function, since the circulation model $\modCirc$ is kept in its full-order form, while only the computationally demanding part, i.e. the 3D electromechanical model $\modEMfom$, is surrogated.
This has a number of advantages:
\begin{enumerate}
    \item The ROM $\modEMred$ is independent of the circulation model $\modCirc$ to which is coupled and it can also be coupled to models different from those used during training, or with different values of its parameters.
    \item Unlike the emulators in \cite{di2018gaussian,dabiri2019prediction,Longobardi2020,Cai2021}, our ROM allows for time extrapolation, i.e. reliable predictions even over longer time spans than those used during training, as demonstrated by our numerical results.
    This observation is very important since, while for emulators that fit the parameters-to-QoIs map the simulations contained in the training set must have reached a limit cycle (otherwise the associated QoIs would be meaningless), with our approach it is possible to defer the reaching of the limit cycle to the ROM, during the online phase.
    This permits to lighten the computational cost associated with the creation of the training set, which typically represents the main component of the total cost (see use cases cited above).
    \item The output of the simulations obtained with our ROM is not limited to a list of QoIs.
    Indeed, it contains the transient of pressures and volumes of the surrogate heart chamber and of the compartments of the circulation model.
    We notice that the latter is not hidden as a black-box in the parameters-to-QoIs map, as for emulators, but is rather represented explicitly.
    \item With our approach, the ROM only needs to learn the variability with respect to $\paramM$ since the circulation model (involving the parameters $\paramC$) remains explicitly represented.
    Conversely, emulators based on a parameters-to-QoIs map must capture the dependence on both $\paramM$ and $\paramC$ and thus the training set must be large enough to accurately represent their statistical variability.
    Indeed, we obtain accurate results with only 30--40 samples, a very low number compared to the ones typically required to construct emulators (e.g. 825 samples in \cite{Longobardi2020}, 9000 samples in \cite{Cai2021}).
\end{enumerate}

A limitation of our ANN-based ROM is that the online phase may be slower than the one of emulators, which do not require to approximate a differential equation but only need to evaluate a function.
However, our ANN-based ROM enables real-time simulations and can be readily applied in a number of use cases, such as global sensitivity analysis, parameter estimation and uncertainty quantification.
Furthermore, as highlighted above, most of the computational cost is not due to the evaluation of the ROM, but rather to the construction of the training set, which our approach allows to keep very small (30--40 samples in our test cases).
Therefore, we conclude that keeping the circulation model in its full-order form and surrogating only the computationally intensive electromechanical model allows for a very favorable trade-off between what is reduced (variability that must be explored during the offline phase) and what is not reduced (variability that must be account for in the online phase).

A different type of emulator is the one we proposed in \cite{regazzoni2021cardioemulator}, that permits - similarly to what is done in this paper - to couple a reduced version of the 3D electromechanical model to a circulation model, thus enabling for real-time numerical simulations.
However, the emulator of \cite{regazzoni2021cardioemulator} is built for a fixed value of the parameters $\paramM$, or (through a linear interpolation) for the variability of a single parameter at most.
On the other hand, its construction only requires one or two numerical simulations performed through the FOM.
Therefore, the emulator of \cite{regazzoni2021cardioemulator} is advantageous when one needs to surrogate the model for a given parametrization (e.g. to quickly converge to a limit cycle, or to perform sensitivity analysis on $\paramC$ only); conversely, when one needs to explore the parametric dependence of an electromechanical model, the approach proposed in this paper turns out to be favorable.

The ROM proposed in this paper features several differences with respect to projection-based ROMs, which have been used in the cardiovascular field as well (see e.g. \cite{QMN16,bonomi2017matrix,pagani2017reduced,pagani2021enabling}), or deep learning based approximations of the parameters-to-solution map (see e.g. \cite{fresca2020deep,fresca2021}).
The latter families of ROMs indeed provide a richer output than our ROM, as they allow for the approximation of spatially varying outputs of the 3D model, such as transmembrane potential and tissue displacement.
However, the online phase of projection-based ROMs is typically more computationally demanding than the one of our ROM \cite{bonomi2017matrix,pagani2021enabling}.
Moreover, unlike our proposed method, projection-based ROMs are intrusive in the equations to be reduced and suitable hyper-reduction techniques should be addressed to handle the nonlinearities of cardiac models. 
The deep learning based method of \cite{fresca2020deep,fresca2021} features a training phase that is more computationally demanding than the one of the method proposed in this paper, due to the larger size of the outputs.
Conversely, in applications in which the knowledge of pressures, volumes and blood fluxes associated with the cardiac chambers are sufficient, our ROM permits to accurately approximate the outputs of the FOM at a very reduced computational cost and in a non-intrusive manner (only pressure and volume recordings are required from the FOM).

We remark that the proposed method is limited to \MS{electromechanical} models that feature a periodic behavior.
The ROM is indeed periodic by construction, due to the presence on the sine and cosine terms, as highlighted in Eq.~\eqref{eqn:model_ANN}.
Therefore, the ROM is not suitable, e.g., to simulate an irregular electrophysiological behavior, such as ventricular arrhythmias, and other electric dysfunctions.

%% file: parts_conclusions.tex
\section{Conclusions}
\label{sec:conclusions}

We presented a machine learning method to build ANN-based ROMs of cardiac \MS{electromechanical} models.
Indeed, on the basis of pressure and volume transients generated with the FOM, we are capable of learning a system of differential equations that approximates the dynamics of the cardiac chamber to be surrogated.
This system of differential equations, linking pressure and volume of a cardiac chamber, is coupled with lumped-parameter models of cardiac hemodynamics, thus allowing for the simulation of the cardiac function at a dramatically reduced computational cost with respect to the original FOM.
As a matter of fact, our ROM permits to perform numerical simulations virtually in real-time.
Moreover, thanks to its non-intrusive nature, the proposed algorithm can be easily applied to other \MS{electromechanical} models besides the one considered in this work.

We presented two \MS{test cases} in which we employ the ANN-based ROM.
We carried out a global sensitivity analysis to assess the influence of the parameters of the electromechanical and hemodynamic models on a list of outputs of clinical interest.
Then, we performed a Bayesian estimation of a couple of parameters (belonging to the electromechanical and hemodynamic models, respectively), starting from the noisy measurement of a couple of scalar quantities (namely maximum and minimum arterial pressure).
In both the cases, performing through the FOM the large number of numerical simulations needed would not have been feasible, due to their high computational cost (it would in fact have taken tens of years on a supercomputer computational platform).
Replacing the FOM with its ANN-based surrogate allowed us to obtain an approximate solution in a few hours of computation.
By taking into account that the generation of the numerical simulations contained in the training set required less than 7 days on the same computational platform, our ANN-based ROM allowed us to reduce the total computational time by more than 3'000 times.

%% file: acknowledgements.tex
\section*{Acknowledgements}
This project has received funding from the European Research Council (ERC) under the European Union's Horizon 2020 research and innovation programme (grant agreement No 740132, iHEART - An Integrated Heart Model for the simulation of the cardiac function, P.I. Prof. A. Quarteroni).
\begin{center}
	\raisebox{-.5\height}{\includegraphics[width=.15\textwidth]{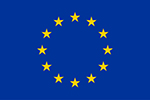}}
	\hspace{2cm}
	\raisebox{-.5\height}{\includegraphics[width=.15\textwidth]{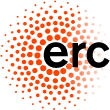}}
\end{center}

%% file: parts_app_params.tex
\section{Numerical settings}
\label{app:params}



In this appendix we report the main numerical \MS{approaches that we employ to develop and test the ROM of cardiac electromechanics presented in this paper.}

\subsection{\MS{Electromechanical} model} \label{app:params:EM}

We use the same \MS{choices} in terms of space and time discretization for all the numerical simulations of this paper.
Specifically, we employ a finer mesh for cardiac electrophysiology, which consists of 258'415 DOFs and 240'864 elements (with an average mesh size of $h_{\text{mean}} = 1.7$ $\si{\milli\meter}$), and a coarser one for cardiac mechanics, which is made by 35'725 DOFs and 30'108 elements ($h_{\text{mean}} = 3.4$ $\si{\milli\meter}$).
To advance the electrophysiological variables we use a timestep length of $\Delta t = \SI{100}{\micro\second}$ and a five times larger timestep to advance the mechanical variables.
The numerical simulations were run by using one cluster node endowed with 32 cores (four Intel Xeon E5-4610 v2, 2.3 GHz) which is available at MOX, Dipartimento di Matematica.

\subsection{Convergence to the limit cycle} \label{app:params:limit-cycle}

To determine when a simulation reached a limit cycle (i.e. a periodic solution), we employ a criterion based on the increment between successive cycles.
Specifically, we consider the limit cycle to be reached when maximum difference between two consecutive cycles in the pressures and volumes of all four chambers is less than \SI{0.8}{\mmHg} and \SI{0.8}{\milli\liter}, respectively.
This criterion is typically satisfied in 5 to 15 cycles.
In any case, we always perform a minimum of 5 cycles.

\subsection{Training algorithm} \label{app:params:ANN}

To generate the training datasets, we employ a Monte Carlo based sampling of the parameter space $\paramSpaceM \times \paramSpaceC$.
As mentioned, we consider a subset of the parameters $\paramC$, made by
\begin{equation*}
    \Vheart, \EpLA, \tLAV, \EpRA, \tRAV, \EaRV, \EpRV, \RarSYS, \RvnSYS.
\end{equation*}
To define the sampling spaces $\paramSpaceM$ and $\paramSpaceC$, we consider 50\% to 200\% ranges of the baseline values reported in Tabs.~\ref{tab:paramsM} and \ref{tab:paramsC} for $\Vheart$, $\EpLA$, $\aXB$, $\EpRA$, $\EaRV$, $\EpRV$, $\RarSYS$ and $\RvnSYS$, 50\% to 150\% for $\EPcondf$ and $\MstiffPass$.
We vary $\fibersAngle \in (\SI{40}{\degree}, \SI{80}{\degree})$ and $\tLAV, \tRAV \in (0.08, 0.24)$ $\si{\second}$.


To train the ANN-based model according to Eq.~\eqref{eqn:optimization_problem}, we employ as mentioned the algorithm that we proposed in \cite{regazzoni2019modellearning}.
More precisely, we approximate the trained ODE with a time step size of $\Delta t = \SI{5}{\milli\second}$ and the objective functional with a twice as wide time step.
To train the ANN, we employ 2000 iterations of the Levenberg-Marquardt algorithm.

\subsection{Hyperparameters tuning} \label{app:params:Xvalidation}

To tune the hyperparameters, we adopt a $k$-fold cross-validation procedure as described in Sec.~\ref{sec:methods:hyperparameters}.
Specifically, we consider a given hyperparameters setting and we partition the training dataset in $k=5$ non-overlapping subsets.
Then, by retaining a single subset as validation data, we train, on the remaining $k-1$ subsets, three different models (with three different random initializations of the ANN weights) and we keep the model attaining the lowest validation error.
This procedure is repeated $k$ times, once for each subset.
Finally, we average the errors over the $k$ trained models, and we compare the resulting average errors for different hyperparameter settings.
More precisely, to compare the performances of the models obtained with different hyperparameter settings, we consider the validation errors associated with the QoIs reported in Tab.~\ref{tab:QoI}, in terms of relative errors.
Moreover, we assess the generalization skills of the trained models by computing the ratio between validation errors and training errors.
Indeed, a ratio that is much larger than one indicates overfitting.
By following these criteria, we get to the final configurations of Tab.~\ref{tab:ROMs}.

\subsection{Global sensitivity analysis} \label{app:params:GSA}

To perform global sensitivity analysis, we consider all parameters $\paramM$ and $\paramC$, which are reported in Tabs.~\ref{tab:paramsM} and~\ref{tab:paramsC}, respectively.
To define the sampling spaces $\paramSpaceM$ and $\paramSpaceC$, we consider 80\% to 120\% ranges of the baseline values in Tabs.~\ref{tab:paramsM} and~\ref{tab:paramsC}, except for $\Vheart \in (193, 593)$ $\si{\milli\liter}$ and for $\tLAV, \tRAV \in (0.04, 0.28)$ $\si{\second}$.

We employ the Saltelli's method to sample the parametric space $\paramSpaceM \times \paramSpaceC$.
A naive sampling of the space would require an exponential increase with respect to the number pf parameters to guarantee a prescribed accuracy; conversely, Saltelli's method leads to a linear increase.
More precisely, the number of samples required is $N (2D + 2)$, where $D$ is the number of parameters and $N$ is a user defined setting.
In this paper we set $N = 1000$, for a total of 74'000 samples, that guarantees small confidence intervals around the first-order Sobol indices and the total-effect Sobol indices reported in Figs.~\ref{fig:SA_all_S1} and~\ref{fig:SA_all_ST}, respectively.
We use one cluster node endowed with 20 cores (four Intel Xeon E5-2640 v4, 2.4 GHz), which is available at MOX, Dipartimento di Matematica, to run global sensitivity analysis.


\subsection{Bayesian parameter estimation} \label{app:params:bayesianPE}

We perform Bayesian parameter estimation by means of the MCMC method.
We consider 500 samples per chain, a jump period of 10 samples and a burn-in equal to 1000, for a total number of 20 chains.
Indeed, we use one cluster node endowed with 20 cores (four Intel Xeon E5-2640 v4, 2.4 GHz), which is available at MOX, Dipartimento di Matematica, where each core manages a single chain and then all final results are collected together at the end of the simulation.